# STATISTICAL INFERENCE FOR TIME-VARYING ARCH PROCESSES[1]

### By Rainer Dahlhaus and Suhasini Subba Rao

*Universität Heidelberg*

In this paper the class of ARCH($\infty$) models is generalized to the nonstationary class of ARCH($\infty$) models with time-varying coefficients. For fixed time points, a stationary approximation is given leading to the notation "locally stationary ARCH($\infty$) process." The asymptotic properties of weighted quasi-likelihood estimators of time-varying ARCH($p$) processes ($p < \infty$) are studied, including asymptotic normality. In particular, the extra bias due to nonstationarity of the process is investigated. Moreover, a Taylor expansion of the nonstationary ARCH process in terms of stationary processes is given and it is proved that the time-varying ARCH process can be written as a time-varying Volterra series.

**1. Introduction.** To model volatility in time series, Engle [6] introduced the ARCH model where the conditional variance is stochastic and dependent on past observations. The ARCH model and several of its related models have gained widespread recognition because they model quite well the volatility in financial markets over relatively short periods of time (cf. [3, 13]). However, underlying all these models is the assumption of stationarity. Now given the changing pace of the world's economy, modeling financial returns over long intervals using stationary time series models may be inappropriate. It is quite plausible that structural changes in financial time series may occur, causing the time series over long intervals to deviate significantly from stationarity. It is therefore plausible that, by relaxing the assumption of stationarity in an adequate way, we may obtain a better fit. In this direction, Drees and Stărică [5] have proposed the simple nonlinear model $X_t = \mu + \sigma(t)Z_t$, where $Z_t$ are independent, identically distributed random variables and $\sigma(\cdot)$ is a smooth function, which they estimate using

Received March 2003; revised July 2005.

[1]Supported in part by the Deutsche Forschungsgemeinschaft (DA 187/12-1).

a nonparametric regression method. Essentially, though it is not mentioned, the authors are treating $\sigma(t)$ as if it were of the form $\sigma(t) = \tilde{\sigma}(t/N)$, with $N$ being the sample size. Through this rescaling device it is possible to obtain a framework for a meaningful asymptotic theory. Feng [7] has also studied time inhomogeneous stochastic volatility, by introducing a multiplicative seasonal and trend component into the GARCH model.

In this paper we generalize the class of ARCH($\infty$) models (cf. [9, 16]) to models with time-varying parameters:

$$(1) \qquad X_t = \sigma(t) Z_t \qquad \text{where } \sigma(t)^2 = a_0(t) + \sum_{j=1}^{\infty} a_j(t) X_{t-j}^2,$$

and $Z_t$ are independent, identically distributed random variables with $\mathbb{E} Z_t = 0$, $\mathbb{E} Z_t^2 = 1$. As in nonparametric regression and in other work on nonparametric statistics, we use the rescaling device to develop an asymptotic theory around such a class of models, that is, we rescale the parameters to the unit interval [see (2) below]. The resulting process is called the time-varying ARCH (tvARCH) process. The same rescaling device has been used, for example, in nonparametric time series by Robinson [15] and by Dahlhaus [4] in his definition of local stationarity which was essentially restricted to time-varying linear processes. We shall show in Section 2 that the tvARCH process can be locally approximated by stationary ARCH processes. Therefore, this new class of tvARCH processes can also be called locally stationary. The stationary ARCH approximation will later be used to transfer results for stationary ARCH processes to the locally stationary situation.

In Section 3 we study parameter estimation for tvARCH($p$) models by weighted quasi-maximum likelihood methods. The nonstationarity of the process causes the estimator to be biased. We will show that the bias can be explained in terms of the derivatives of the tvARCH process. Furthermore, we will prove asymptotic normality of the estimator.

In Section 4 we also define a special derivative of the tvARCH process and give a Taylor expansion of the nonstationary tvARCH process in terms of stationary processes. This derivative enables us to study more precisely the nonstationary behavior of the process. Moreover, the derivative process turns out to be a solution of a stochastic differential equation.

In Section 5 time-varying Volterra series are studied. They are used to prove the existence of a tvARCH($\infty$) process and to derive the results of Section 4 on its derivatives. It is worth noting that the results in Section 5 are of independent interest and the methods used here can be generalized to other nonstationary processes.

In the Appendix we prove convergence theorems for ergodic stationary processes and some specific convergence and approximation results for the likelihood process. We also derive mixing properties of several processes, including derivatives of the likelihood process.

**2. The time varying ARCH process.** In this section we broaden the class of ARCH($\infty$) models, by introducing nonstationary ARCH($\infty$) models with time-dependent parameters. In order to obtain a framework for a meaningful asymptotic theory, we rescale the parameter functions as in nonparametric regression and for (linear) locally stationary processes to the unit interval, that is, we assume

(2)
$$X_{t,N} = \sigma_{t,N} Z_t$$
$$\text{where } \sigma_{t,N}^2 = a_0\left(\frac{t}{N}\right) + \sum_{j=1}^{\infty} a_j\left(\frac{t}{N}\right) X_{t-j,N}^2 \text{ for } t = 1, \ldots, N,$$

where $Z_t$ are independent, identically distributed random variables with $\mathbb{E} Z_t = 0$, $\mathbb{E} Z_t^2 = 1$. We call the sequence of stochastic processes $\{X_{t,N} : t = 1, \ldots, N\}$ which satisfy (2) a time-varying ARCH (tvARCH) process. As shown below, the tvARCH-process can be locally approximated by stationary ARCH processes. Therefore, we also call tvARCH processes locally stationary.

We mention that the rescaling technique is mainly introduced for obtaining a meaningful asymptotic theory, and by this device we can obtain adequate approximations for the nonrescaled case. In particular, the rescaling does not effect the estimation procedure. Furthermore, classical ARCH models are included as a special case (if the parameters are constant in time).

We make the following assumptions.

ASSUMPTION 1. The sequence of stochastic processes $\{X_{t,N} : t = 1, \ldots, N\}$ has a time-varying ARCH representation defined in (2) where the parameters satisfy the following properties: There exist constants $0 < \rho, Q, M < \infty$, $0 < \nu < 1$ and a positive sequence $\{\ell(j)\}$ such that $\inf_u a_0(u) > \rho$ and

(3)
$$\sup_u a_j(u) \leq \frac{Q}{\ell(j)},$$

(4)
$$Q \sum_{j=1}^{\infty} \frac{1}{\ell(j)} \leq (1 - \nu),$$

(5)
$$|a_j(u) - a_j(v)| \leq M \frac{|u - v|}{\ell(j)},$$

where $\{\ell(j)\}$ satisfies

$$\sum_{j \geq 1} \frac{j}{\ell(j)} < \infty.$$

An example of such a positive sequence $\{\ell(j)\}$ is

$$\ell(j) = \begin{cases} 1, & j = 1, \\ j^2 \log^{1+\kappa} j, & j > 1, \end{cases}$$

with some $\kappa > 0$ or $\ell(j) = \eta^j$ for some $\eta > 1$. Condition (4) implies that $\mathbb{E}(X_{t,N}^2)$ is uniformly bounded over $t$ and $N$.

PROPOSITION 1. *Under Assumption* 1, $\{X_{t,N}^2\}$ *defined in* (2) *has an almost surely well-defined unique solution in the set of all causal solutions. The solution has the form of a time-varying Volterra series expansion.*

The proof for Proposition 1, as well as all the other proofs of results in this section can be found in Section 5. We mention that a similar result also holds for nonrescaled tvARCH($\infty$) processes.

It is worth noting that throughout this paper we shall be working with $X_{t,N}^2$ rather than $X_{t,N}$, unless stated otherwise. This is because $X_{t,N}$ can be randomly either positive or negative, whereas $X_{t,N}^2$ is always positive, allowing it to be unique.

The smoothness of the parameters $\{a_j(\cdot)\}$ guarantees that the process has (asymptotically) locally a stationary behavior. We now make this notion precise. The first point of interest is to study the stationary process which locally approximates the tvARCH-process in some neighborhood of a fixed point $t_0$ (or in rescaled time $u_0$). For each given $u_0 \in (0, 1]$, the stochastic process $\{\tilde{X}_t(u_0)\}$ is the stationary ARCH process associated with the tvARCH($\infty$) process at time point $u_0$ if it satisfies

(6)
$$\tilde{X}_t(u_0) = \sigma_t(u_0) Z_t,$$

$$\text{where } \sigma_t(u_0)^2 = a_0(u_0) + \sum_{j=1}^{\infty} a_j(u_0) \tilde{X}_{t-j}(u_0)^2$$

for all $t \in \mathbb{Z}$. It is worth noting, if the parameters $\{a_j(u_0)\}$ satisfy Assumption 1, then $\{\tilde{X}_t(u_0)\}$ is a stationary, ergodic ARCH($\infty$) process (cf. [9]).

Comparing (6) with (2), it seems clear that if $t/N$ is close to $u_0$, then $X_{t,N}^2$ and $\tilde{X}_t(u_0)^2$ should be close and the degree of the approximation should depend both on the rescaling factor $N$ and the deviation $|t/N - u_0|$. This is shown below.

THEOREM 1. *Suppose* $\{X_{t,N}\}$ *is a tvARCH process which satisfies Assumption* 1 *and let* $\tilde{X}_t(u_0)$ *be defined as in* (6). *Then there exist a stationary, ergodic, positive process* $\{U_t\}$ *independent of* $u_0$ *with finite mean and a constant* $K$ *independent of* $t$ *and* $N$ *such that*

(7)
$$|X_{t,N}^2 - \tilde{X}_t(u_0)^2| \leq K \left( \left| \frac{t}{N} - u_0 \right| + \frac{1}{N} \right) U_t \qquad a.s.$$

We mention that an explicit formula for $U_t$ is given in (48). As a consequence of (7), we have

$$X_{t,N}^2 = \tilde{X}_t(u_0)^2 + O_p\left(\left|\frac{t}{N} - u_0\right| + \frac{1}{N}\right).$$

The bound in (7) allows us to approximate the local average of $X_{t,N}^2$ by an average of $\tilde{X}_t(u_0)^2$ (this is of particular interest here, since the local average and weighted local average will be used frequently in later sections). For example, suppose $|u_0 - t_0/N| < 1/N$ and we average $X_{t,N}^2$ about a neighborhood whose length $(2M + 1)$ increases as $N$ increases but where the ratio $M/N \to 0$ as $N \to \infty$. Then by using (7), we have

$$(8) \qquad \frac{1}{2M+1}\sum_{k=-M}^{M} X_{t_0+k,N}^2 = \frac{1}{2M+1}\sum_{k=-M}^{M} \tilde{X}_{t_0+k}(u_0)^2 + R_N,$$

where $R_N$ is bounded by

$$|R_N| \leq K\frac{M}{N}\left\{\frac{1}{2M+1}\sum_{k=-M}^{M} U_{t_0+k}\right\} \xrightarrow{\mathcal{P}} 0.$$

Thus, about the time point $t_0$ the local average of a tvARCH process is asymptotically the same as the local average of the stationary ARCH process $\{\tilde{X}_t(u_0)^2\}$. Therefore, by using (7), we can locally approximate the tvARCH process by a stationary process. The above approximation can be refined by using derivative processes as defined in Section 4. By using them, we can find, for example, an expression for the asymptotic bias $R_N$ in (8).

**3. The segment quasi-likelihood estimate.** In this section we consider a kernel type estimator of the parameters of a tvARCH($p$) model given the sample $\{X_{t,N} : t = 1, \ldots, N\}$. The process $\{X_{t,N}\}$ is assumed to satisfy the representation

$$X_{t,N} = \sigma_{t,N} Z_t,$$

$$(9) \qquad \text{where } \sigma_{t,N}^2 = a_0\left(\frac{t}{N}\right) + \sum_{j=1}^{p} a_j\left(\frac{t}{N}\right) X_{t-j,N}^2 \text{ for } t = 1 \ldots, N,$$

where $Z_t$ are independent, identically distributed random variables with $\mathbb{E}Z_t = 0$, $\mathbb{E}Z_t^2 = 1$. The order $p$ is assumed known. We study the distributional properties of the estimator, including asymptotic normality. Furthermore, we will investigate the bias of the estimator due to nonstationarity of the tvARCH($p$) process. We will use the following assumptions.

ASSUMPTION 2. The sequence of stochastic processes $\{X_{t,N} : t = 1, \ldots, N\}$ has a tvARCH($p$) representation defined by (9). Furthermore:

(i) The process satisfies Assumption 1.

(ii) For some $\delta > 0$,

$$\mathbb{E}(|Z_t|^{4(1+\delta)}) < \infty. \tag{10}$$

(iii) Let $\Omega$ be the compact set

$$\Omega = \left\{ \boldsymbol{\alpha} = (\alpha_0, \alpha_1, \dots, \alpha_p) : \sum_{j=1}^{p} \alpha_j \leq 1, \rho_1 \leq \alpha_0 \leq \rho_2, \rho_1 \leq \alpha_i \text{ for } i = 1, \dots, p \right\},$$

where $0 < \rho_1 \leq \rho_2 < \infty$. For each $u \in (0, 1]$, we assume $\mathbf{a}_u \in \text{Int}(\Omega)$, where $\mathbf{a}_u = (a_0(u), a_1(u), \dots, a_p(u))$.

(iv) The third derivative of $a_j(\cdot)$ exists with

$$\sup_u \left| \frac{\partial^i a_j(u)}{\partial u^i} \right| \leq C$$

for $i = 1, 2, 3$ and $j = 0, 1, \dots, p$, where $C$ is a finite constant independent of $i$ and $j$.

(v) The random variable $Z_t$ has a positive density on an interval containing zero.

(vi) [This assumption is only used in Theorem 3(ii)]

$$\{\mathbb{E}(Z_0^{12})\}^{1/6} \sum_{j=1}^{p} \frac{Q}{\ell(j)} \leq (1 - \nu).$$

REMARK 1.  (i) The conditions placed on the parameter space in Assumption 2(iii) can be relaxed to include all vectors $\boldsymbol{\alpha} = (\alpha_0, \alpha_1, \dots, \alpha_p)$, where $\alpha_i = 0$ for any $i = 1, \dots, p$, in the parameter space. By including these points, a method for model selection could be derived. However, the cost for relaxing this assumption is that additional moment conditions have to be placed on $X_{t,N}$.

(ii) We use Assumption 2(ii) to prove asymptotic normality of the estimator. Typically for stationary ARCH processes, the result can be proved if $\mathbb{E}(Z_t^4) < \infty$. However, we require the mildly stronger assumption $\mathbb{E}(|Z_t^{4+\delta}) < \infty$ to prove a similar result for sums of martingale arrays as opposed to sums of martingale differences used in the stationary situation (cf. [10], Theorem 3.2). Assumption 2(vi) means that both $\mathbb{E}(X_{t,N}^{12})$ and $\mathbb{E}(\tilde{X}_t(u)^{12})$ are uniformly bounded in $t, N$ and $u$. We refer also to the comments on the moment assumptions in Section 6.

(iii) In Section 5 we apply a theorem of Basrak, Davis and Mikosch [1], who gave conditions under which a GARCH$(p, q)$ process is mixing. Assumption 2(iii) is sufficient for the Lyapunov exponent of the random recurrence matrix associated with $\{\tilde{X}_t(u)\}$ to be negative. In addition, Assumption 2(iii), (v) is sufficient to ensure that the ARCH process $\{\tilde{X}_t(u)\}$ is $\alpha$-mixing of rate $-\infty$ (see [1] and references therein).

We now define the segment (kernel) estimator of $\boldsymbol{a}(u_0)$ for each $u_0 \in (0,1)$. Let $t_0 \in \mathbb{N}$ such that $|u_0 - t_0/N| < 1/N$. The estimator considered in this section is the minimizer of the weighted conditional likelihood

$$\mathcal{L}_{t_0, N}(\boldsymbol{\alpha}) := \sum_{k=p+1}^{N} \frac{1}{bN} W\left(\frac{t_0 - k}{bN}\right) \ell_{k,N}(\boldsymbol{\alpha}), \tag{11}$$

where

$$\ell_{k,N}(\boldsymbol{\alpha}) = \frac{1}{2}\left(\log w_{k,N}(\boldsymbol{\alpha}) + \frac{X_{k,N}^2}{w_{k,N}(\boldsymbol{\alpha})}\right) \tag{12}$$

$$\text{with } w_{k,N}(\boldsymbol{\alpha}) = \alpha_0 + \sum_{j=1}^{p} \alpha_j X_{k-j,N}^2$$

and $W : [-1/2, 1/2] \to \mathbb{R}$ is a kernel function of bounded variation with $\int_{-1/2}^{1/2} W(x)\,dx = 1$ and $\int_{-1/2}^{1/2} x W(x)\,dx = 0$. That is, we consider

$$\hat{\mathbf{a}}_{t_0, N} = \underset{\boldsymbol{\alpha} \in \Omega}{\arg\min} \, \mathcal{L}_{t_0, N}(\boldsymbol{\alpha}). \tag{13}$$

Obviously $\ell_{t,N}(\boldsymbol{\alpha})$ is the conditional likelihood of $X_{t,N}$ given $X_{t-1,N}, \ldots,$ $X_{t-p,N}$ and the parameters $\boldsymbol{\alpha} = (\alpha_0, \ldots, \alpha_p)^T$, provided the $Z_t$ are normally distributed. All results below also hold if the $Z_t$ are not normally distributed but simply satisfy Assumption 2. For this reason (and the fact that the conditional likelihood is not the full likelihood), the likelihood is called a quasi-likelihood. For later reference, we list the derivatives of $\ell_{k,N}(\boldsymbol{\alpha})$. Let $\nabla = (\frac{\partial}{\partial \alpha_0}, \ldots, \frac{\partial}{\partial \alpha_p})^T$. Since $\nabla^2 w_{k,N}(\boldsymbol{\alpha}) = 0$, we have

$$\ell_{k,N}(\boldsymbol{\alpha}) = \frac{1}{2}\left\{\log(w_{k,N}(\boldsymbol{\alpha})) + \frac{X_{k,N}^2}{w_{k,N}(\boldsymbol{\alpha})}\right\}, \tag{14}$$

$$\nabla \ell_{k,N}(\boldsymbol{\alpha}) = \frac{1}{2}\left\{\frac{\nabla w_{k,N}(\boldsymbol{\alpha})}{w_{k,N}(\boldsymbol{\alpha})} - \frac{X_{k,N}^2 \nabla w_{k,N}(\boldsymbol{\alpha})}{w_{k,N}(\boldsymbol{\alpha})^2}\right\}, \tag{15}$$

$$\begin{aligned} \nabla^2 \ell_{k,N}(\boldsymbol{\alpha}) = \frac{1}{2}\bigg\{ &-\frac{\nabla w_{k,N}(\boldsymbol{\alpha}) \nabla w_{k,N}(\boldsymbol{\alpha})^T}{w_{k,N}(\boldsymbol{\alpha})^2} \\ &+ 2\frac{X_{k,N}^2 \nabla w_{k,N}(\boldsymbol{\alpha}) \nabla w_{k,N}(\boldsymbol{\alpha})^T}{w_{k,N}(\boldsymbol{\alpha})^3}\bigg\}. \end{aligned} \tag{16}$$

$\hat{\mathbf{a}}_{t_0, N}$ is regarded as an estimator of $\mathbf{a}_{t_0/N} = (a_0(t_0/N), \ldots, a_p(t_0/N))^T$ or of $\mathbf{a}_{u_0}$, where $|u_0 - t_0/N| < 1/N$.

In the derivation of the asymptotic properties of this estimator we make use of the local approximation of $X_{t,N}^2$ by the stationary process $\tilde{X}_t(u_0)^2$ defined in Section 2. Similarly to the above, we therefore define the weighted likelihood

$$\tilde{\mathcal{L}}_N(u_0, \boldsymbol{\alpha}) := \sum_{k=p+1}^{N} \frac{1}{bN} W\left(\frac{t_0 - k}{bN}\right) \tilde{\ell}_k(u_0, \boldsymbol{\alpha}), \tag{17}$$

where $|u_0 - t_0/N| < 1/N$ and

$$\tilde{\ell}_t(u_0, \boldsymbol{\alpha}) = \frac{1}{2} \left( \log \tilde{w}_t(u_0, \boldsymbol{\alpha}) + \frac{\tilde{X}_t(u_0)^2}{\tilde{w}_t(u_0, \boldsymbol{\alpha})} \right) \tag{18}$$

$$\text{with } \tilde{w}_t(u_0, \boldsymbol{\alpha}) = \alpha_0 + \sum_{j=1}^{p} \alpha_j \tilde{X}_{t-j}(u_0)^2.$$

It is obvious that the same formulas as in (14)–(16) also hold for $\tilde{\ell}_k(u_0, \boldsymbol{\alpha})$ with $X_{k,N}^2$ and $w_{k,N}(\boldsymbol{\alpha})$ replaced by $\tilde{X}_k(u_0)^2$ and $\tilde{w}_k(u_0, \boldsymbol{\alpha})$, respectively.

It is shown below that both $\mathcal{L}_{t_0,N}(\boldsymbol{\alpha})$ and $\tilde{\mathcal{L}}_N(u_0, \boldsymbol{\alpha})$ converge to

$$\mathcal{L}(u_0, \boldsymbol{\alpha}) := \mathbb{E}(\tilde{\ell}_0(u_0, \boldsymbol{\alpha})) \tag{19}$$

as $N \to \infty$, $b \to 0$, $bN \to \infty$ and $|u_0 - t_0/N| < 1/N$. It is easy to show that $\mathcal{L}(u_0, \boldsymbol{\alpha})$ is minimized by $\boldsymbol{\alpha} = \mathbf{a}_{u_0}$.

Furthermore, let

$$\mathcal{B}_{t_0,N}(\boldsymbol{\alpha}) := \mathcal{L}_{t_0,N}(\boldsymbol{\alpha}) - \tilde{\mathcal{L}}_N(u_0, \boldsymbol{\alpha})$$
$$= \sum_{k=p+1}^{N} \frac{1}{bN} W\left(\frac{t_0 - k}{bN}\right) (\ell_{k,N}(\boldsymbol{\alpha}) - \tilde{\ell}_k(u_0, \boldsymbol{\alpha})). \tag{20}$$

Since $\tilde{\mathcal{L}}_N(u_0, \boldsymbol{\alpha})$ is the likelihood of the stationary approximation, $\tilde{X}_t(u_0)\mathcal{B}_{t_0,N}(\boldsymbol{\alpha})$ is a bias caused by the deviation from stationarity. Lemma A.6 implies that $\mathcal{B}_{t_0,N}(\boldsymbol{\alpha}) = O_p(b)$. A better rate will be derived by a Taylor expansion in Proposition 3. Let

$$\Sigma(u_0) = \frac{1}{2} \mathbb{E}\left\{ \frac{\nabla \tilde{w}_0(u_0, \mathbf{a}_{u_0}) \nabla \tilde{w}_0(u_0, \mathbf{a}_{u_0})^T}{\tilde{w}_0(u_0, \mathbf{a}_{u_0})^2} \right\}. \tag{21}$$

Since $\tilde{X}_k(u_0)/\tilde{w}_k(u_0, \mathbf{a}_0) = Z_k^2$ and $Z_k^2$ is independent of $\tilde{w}_k(u_0, \mathbf{a}_0)$, we have

$$\mathbb{E}(\nabla^2 \tilde{\ell}_0(u_0, \mathbf{a}_{u_0})) = -\Sigma(u_0)$$

and

$$\mathbb{E}(\nabla \tilde{\ell}_0(u_0, \mathbf{a}_{u_0}) \nabla \tilde{\ell}_0(u_0, \mathbf{a}_{u_0})^T) = \frac{\text{var}(Z_0^2)}{2} \Sigma(u_0).$$

If $Z_t$ is Gaussian, then $\text{var}(Z_0^2) = 2$.

Lemma 1. *Suppose $\{X_{t,N} : t = 1, \ldots, N\}$ is a* tvARCH$(p)$ *process which satisfies Assumption* 2(i), (iii) *and let* $\tilde{\mathcal{L}}_N(u_0, \boldsymbol{\alpha})$, $\mathcal{L}(u_0, \boldsymbol{\alpha})$ *and* $\mathcal{B}_{t_0,N}$ *be as defined by* (17), (19) *and* (20), *respectively. Then*

$$\sup_{\boldsymbol{\alpha} \in \Omega} |\tilde{\mathcal{L}}_N(u_0, \boldsymbol{\alpha}) - \mathcal{L}(u_0, \boldsymbol{\alpha})| \xrightarrow{\mathcal{P}} 0, \tag{22}$$

$$\sup_{\boldsymbol{\alpha} \in \Omega} |\mathcal{B}_{t_0,N}(\boldsymbol{\alpha})| \xrightarrow{\mathcal{P}} 0, \tag{23}$$

$$\sup_{\boldsymbol{\alpha} \in \Omega} |\nabla^2 \tilde{\mathcal{L}}_N(u_0, \boldsymbol{\alpha}) - \nabla^2 \mathcal{L}(u_0, \boldsymbol{\alpha})| \xrightarrow{\mathcal{P}} 0 \tag{24}$$

*and*

$$\sup_{\boldsymbol{\alpha} \in \Omega} |\nabla^2 \mathcal{B}_{t_0,N}(\boldsymbol{\alpha})| \xrightarrow{\mathcal{P}} 0, \tag{25}$$

*where $b \to 0$ and $bN \to \infty$ as $N \to \infty$.*

A direct implication of the lemma above is the following corollary.

Corollary 1. *Let $\mathcal{L}_{t,N}(\boldsymbol{\alpha})$ be as defined as in* (11). *Then under the assumptions in Lemma* 1, *we have*

$$\sup_{\boldsymbol{\alpha} \in \Omega} |\mathcal{L}_{t_0,N}(\boldsymbol{\alpha}) - \mathcal{L}(u_0, \boldsymbol{\alpha})| \xrightarrow{\mathcal{P}} 0 \tag{26}$$

*and*

$$\sup_{\boldsymbol{\alpha} \in \Omega} |\nabla^2 \mathcal{L}_{t_0,N}(\boldsymbol{\alpha}) - \nabla^2 \mathcal{L}(u_0, \boldsymbol{\alpha})| \xrightarrow{\mathcal{P}} 0, \tag{27}$$

*where $b \to 0$, $bN \to \infty$ as $N \to \infty$.*

In the theorem below we show that $\hat{\mathbf{a}}_{t_0,N}$ is a consistent estimator of $\mathbf{a}_{u_0}$.

Theorem 2. *Suppose $\{X_{t,N} : t = 1, \ldots, N\}$ is a* tvARCH$(p)$ *process which satisfies Assumption* 2(i), (iii) *and the estimator $\hat{\mathbf{a}}_{t_0,N}$ is as defined in* (13). *Then if $|u_0 - t_0/N| < 1/N$, we have*

$$\hat{\mathbf{a}}_{t_0,N} \xrightarrow{\mathcal{P}} \mathbf{a}(u_0),$$

*where $b \to 0$ and $bN \to \infty$ as $N \to \infty$.*

Proof. By using (26), we have pointwise convergence $\mathcal{L}_{t_0,N}(\mathbf{a}) \xrightarrow{\mathcal{P}} \mathcal{L}(u_0, \mathbf{a})$. Since $\mathbf{a}_{u_0} = \arg\min_{\boldsymbol{\alpha}} \mathcal{L}(u_0, \boldsymbol{\alpha})$, we have

$$\mathcal{L}_{t_0,N}(\hat{\mathbf{a}}_{t_0,N}) \le \mathcal{L}_{t_0,N}(\mathbf{a}_{u_0}) \xrightarrow{\mathcal{P}} \mathcal{L}(u_0, \mathbf{a}_{u_0}) \le \mathcal{L}(u_0, \hat{\mathbf{a}}_{t_0,N}).$$

With (26), we now obtain $\mathcal{L}_{t_0,N}(\hat{\mathbf{a}}_{t_0,N}) \xrightarrow{\mathcal{P}} \mathcal{L}(u_0, \mathbf{a}_{u_0})$. From the continuity of $\mathcal{L}(u_0, \cdot)$ and the compactness of the parameter space, we can now conclude $\hat{\mathbf{a}}_{t_0,N} \xrightarrow{\mathcal{P}} \mathbf{a}_{u_0}$, provided $\mathcal{L}(u_0, \cdot)$ has a unique solution. Since $\mathcal{L}(u_0, \cdot)$ is the same function as in the stationary case, this follows from Lemma 5.5 of [2]. $\square$

We now prove asymptotic normality of the estimator with the usual Taylor expansion argument. We have

$$(28) \quad \nabla\mathcal{L}_{t_0,N}(\hat{\mathbf{a}}_{t_0,N})_i - \nabla\mathcal{L}_{t_0,N}(\mathbf{a}_{u_0})_i = \{\nabla^2\mathcal{L}_{t_0,N}(\bar{\mathbf{a}}_{t_0,N}^i)(\hat{\mathbf{a}}_{t_0,N} - \mathbf{a}_{u_0})\}_i,$$

with $\bar{\mathbf{a}}_{t_0,N}^i$ between $\hat{\mathbf{a}}_{t_0,N}$ and $\mathbf{a}_{u_0}$. Since $\mathbf{a}_{u_0}$ is in the interior of $\Omega$, we have $\sqrt{bN}\nabla\mathcal{L}_{t_0,N}(\hat{\mathbf{a}}_{t_0,N})_i \xrightarrow{\mathcal{P}} 0$. Since $\bar{\mathbf{a}}_{t_0,N}^i \xrightarrow{\mathcal{P}} \mathbf{a}_{u_0}$ and $\sup_{\boldsymbol{\alpha}\in\Omega}|\nabla^2\mathcal{L}_{t_0,N}(\boldsymbol{\alpha}) - \nabla^2\mathcal{L}(u_0, \boldsymbol{\alpha})| \xrightarrow{\mathcal{P}} 0$ [see (27)], then $\nabla^2\mathcal{L}_{t_0,N}(\bar{\mathbf{a}}_{t_0,N}^i) \xrightarrow{\mathcal{P}} -\Sigma(u_0)$. Note that $\Sigma(u_0)$ is nonsingular. This follows from Lemma 5.7 of [2] since $\Sigma(u_0)$ is the same as in the stationary case.

Therefore, the distributional properties of $\hat{\mathbf{a}}_{t_0,N} - \mathbf{a}_{u_0}$ are determined by $\nabla\mathcal{L}_{t_0,N}(\mathbf{a}_{u_0})$. By using (20), we see that

$$(29) \quad \nabla\mathcal{L}_{t_0,N}(\mathbf{a}_{u_0}) = \nabla\tilde{\mathcal{L}}_N(u_0, \mathbf{a}_{u_0}) + \nabla\mathcal{B}_{t_0,N}(\mathbf{a}_{u_0}),$$

which is essentially a decomposition into a stochastic and a bias part [although $\nabla\mathcal{B}_{t_0,N}(\mathbf{a}_{u_0})$ is also random, but its variance is of a lower order—see the details below]. The bias measures the deviation from stationarity and will disappear for a suitable choice of bandwidth $b$ (see Proposition 2 and Theorem 3 below). By substituting (29) into (28), we have

$$\sqrt{bN}((\hat{\mathbf{a}}_{t_0,N} - \mathbf{a}_{u_0}) + \Sigma(u_0)^{-1}\nabla\mathcal{B}_{t_0,N}(\mathbf{a}_{u_0}))$$
$$= \sqrt{bN}\Sigma(u_0)^{-1}\nabla\tilde{\mathcal{L}}_N(u_0, \mathbf{a}_{u_0}) + o_p(1).$$

Thus, the asymptotic distribution of $(\hat{\mathbf{a}}_{t_0,N} - \mathbf{a}_{u_0})$ is determined by $\nabla\tilde{\mathcal{L}}_N(u_0, \mathbf{a}_{u_0})$. Note that this is the gradient of the likelihood of the stationary process $\tilde{X}_t(u_0)^2$ however, with kernel weights. Since

$$(30) \quad \nabla\tilde{\ell}_k(u_0, \mathbf{a}_{u_0}) = \frac{1}{2}\frac{(1 - Z_k^2)\nabla\tilde{w}_k(u_0, \mathbf{a}_{u_0})}{\tilde{w}_k(u_0, \mathbf{a}_{u_0})}$$

is a martingale difference, $\nabla\tilde{\mathcal{L}}_N(u_0, \mathbf{a}_{u_0})$ is the weighted sum of martingale differences.

PROPOSITION 2. *Suppose* $\{X_{t,N} : t = 1, \ldots, N\}$ *is a* tvARCH($p$) *process which satisfies Assumption 2*(i), (ii), (iii) *and* $\tilde{\mathcal{L}}_N(u_0, \mathbf{a}_{u_0})$ *is as defined in* (17). *Then if* $|u_0 - t_0/N| < 1/N$ *we have*

$$(31) \quad \sqrt{bN}\nabla\tilde{\mathcal{L}}_N(u_0, \mathbf{a}_{u_0}) \xrightarrow{\mathcal{D}} \mathcal{N}\left(0, w_2\frac{\mathrm{var}(Z_0^2)}{2}\Sigma(u_0)\right),$$

*where $b \to 0$, $bN \to \infty$, $N \to \infty$ and $w_2 = \int_{-1/2}^{1/2} W(x)^2 \, dx$.*

PROOF. Since $\nabla \tilde{\mathcal{L}}_N(u_0, \mathbf{a}_{u_0})$ is the weighted sum of martingale differences, the result follows from the martingale central limit and the Cramér–Wold device. It is straightforward to check the conditional Lindeberg condition and the conditional variance condition. We omit the details. □

We now consider the stochastic bias $\nabla \mathcal{B}_{t_0, N}(\mathbf{a}_{u_0})$. By using (30) and Lemma A.6, we immediately get the relation

$$\nabla \mathcal{B}_{t_0, N}(\mathbf{a}_{u_0}) = O_p(b). \tag{32}$$

This bound together with the Proposition 2 leads to the assertion of Theorem 3(i) below. As mentioned above, the stochastic bias is a measure for the deviation of the process $\{\ell_{t,N}(\mathbf{a}_{u_0})\}$ from stationarity. This deviation depends on the rate of change of the parameters $\{a_j(u)\}$. Under stronger moment conditions, we will now determine this bias. To achieve this, we replace $\nabla \ell_{k,N}(\mathbf{a}_{u_0})$ by $\nabla \tilde{\ell}_k(\frac{k}{N}, \mathbf{a}_{u_0})$:

$$\nabla \mathcal{B}_{t_0, N}(\mathbf{a}_{u_0}) = \sum_k \frac{1}{bN} W\left(\frac{t_0 - k}{bN}\right)\left(\nabla \tilde{\ell}_k\left(\frac{k}{N}, \mathbf{a}_{u_0}\right) - \nabla \tilde{\ell}_k(u_0, \mathbf{a}_{u_0})\right) + R_N, \tag{33}$$

where

$$R_N = \sum_k \frac{1}{bN} W\left(\frac{t_0 - k}{bN}\right)\left(\nabla \ell_{k,N}(\mathbf{a}_{u_0}) - \nabla \tilde{\ell}_k\left(\frac{k}{N}, \mathbf{a}_{u_0}\right)\right).$$

Corollary A.1 now implies

$$\left|\nabla \ell_{k,N}(\mathbf{a}_{u_0}) - \nabla \tilde{\ell}_k\left(\frac{k}{N}, \mathbf{a}_{u_0}\right)\right| \le \frac{K}{N}\left(U_k + (1 + Z_k^2)\sum_{j=1}^{p} U_{k-j}\right),$$

with some constant $K$ uniformly in $k$. Lemma 1 together with the independence of $Z_k^2$ and $U_{k-j}$ now imply

$$\mathbb{E}(R_N^2)^{1/2} = O\left(\frac{1}{N}\right).$$

Suppose for each $j = 0, \ldots, p$ the third derivative of $a_j(\cdot)$ exists and is uniformly bounded. Then by using Corollary 3 and taking a Taylor expansion of $\nabla \tilde{\ell}_k(u, \mathbf{a}_{u_0})$ about $u = u_0$, we have

$$\begin{aligned}
\nabla \tilde{\ell}_k\left(\frac{k}{N}, \mathbf{a}_{u_0}\right) - \nabla \tilde{\ell}_k(u_0, \mathbf{a}_{u_0}) = {}&\left(\frac{k}{N} - u_0\right)\frac{\partial \nabla \tilde{\ell}_k(u, \mathbf{a}_{u_0})}{\partial u}\bigg|_{u=u_0} \\
&+ \frac{(k/N - u_0)^2}{2}\frac{\partial^2 \nabla \tilde{\ell}_k(u, \mathbf{a}_{u_0})}{\partial u^2}\bigg|_{u=u_0} \\
&+ \frac{(k/N - u_0)^3}{3!}\frac{\partial^3 \nabla \tilde{\ell}_k(u, \mathbf{a}_{u_0})}{\partial u^3}\bigg|_{u=\tilde{U}_k},
\end{aligned} \tag{34}$$

where the random variable $\tilde{U}_k \in (0,1]$. A detailed investigation of the different terms now leads to the following result on $\nabla \mathcal{B}_{t_0,N}(\mathbf{a}_{u_0})$. We mention that, in particular, the expectation of the first term cancels out.

PROPOSITION 3. *Suppose* $\{X_{t,N} : t = 1, \ldots, N\}$ *is a* tvARCH($p$) *process which satisfies Assumption* 2 *and* $W$ *is a kernel function of bounded variation with* $\int_{-1/2}^{1/2} W(x)\,dx = 1$ *and* $\int_{-1/2}^{1/2} W(x)x\,dx = 0$. *Then if* $|u_0 - t_0/N| < 1/N$, *we have*

$$\mathbb{E}(\nabla \mathcal{B}_{t_0,N}(\mathbf{a}_{u_0})) = \frac{1}{2}b^2 w(2)\frac{\partial^2 \nabla \mathcal{L}(u, \mathbf{a}_{u_0})}{\partial u^2}\Big|_{u=u_0} + O\left(b^3 + \frac{1}{N}\right)$$

*and*

$$\mathrm{var}(\nabla \mathcal{B}_{t_0,N}(\mathbf{a}_{u_0})) = O\left(b^6 + \frac{1}{N}\right),$$

*where* $w(2) = \int_{-1/2}^{1/2} W(x)x^2\,dx$.

A detailed proof can be found in Appendix A.4.

Propositions 2 and 3 and (32) give us the distributional properties of the estimator $\hat{\mathbf{a}}_{t_0,N}$, which we summarize in the theorem below.

THEOREM 3. *Suppose* $\{X_{t,N} : t = 1, \ldots, N\}$ *is a* tvARCH($p$) *process which satisfies Assumption* 2(i), (ii), (iii) *and* $W$ *is a kernel function of bounded variation with* $\int_{-1/2}^{1/2} W(x)\,dx = 1$ *and* $\int_{-1/2}^{1/2} W(x)x\,dx = 0$. *Then if* $|u_0 - t_0/N| < 1/N$, *we have the following:*

(i) *If* $b^3 \ll N^{-1}$, *then* $\sqrt{bN}\mathcal{B}_{t_0,N}(\mathbf{a}_{u_0}) \overset{\mathcal{P}}{\to} 0$ *and*

$$\sqrt{bN}(\hat{\mathbf{a}}_{t_0,N} - \mathbf{a}_{u_0}) \overset{\mathcal{D}}{\to} \mathcal{N}\left(0, w_2\frac{\mathrm{var}(Z_0^2)}{2}\Sigma(u_0)^{-1}\right).$$

(ii) *If in addition Assumption* 2(iv), (v), (vi) *holds and* $b^{13} \ll N^{-1}$, *then*

$$\sqrt{bN}\Sigma(u_0)^{-1}\nabla \mathcal{B}_{t_0,N}(\mathbf{a}_{u_0}) = \sqrt{bN}b^2\mu(u_0) + o_p(1)$$

*and*

(35)   $$\sqrt{bN}(\hat{\mathbf{a}}_{t_0,N} - \mathbf{a}_{u_0}) + \sqrt{bN}b^2\mu(u_0) \overset{\mathcal{D}}{\to} \mathcal{N}\left(0, w_2\frac{\mathrm{var}(Z_0^2)}{2}\Sigma(u_0)^{-1}\right),$$

*where*

(36)   $$\mu(u_0) = \frac{1}{2}w(2)\Sigma(u_0)^{-1}\frac{\partial^2 \nabla \mathcal{L}(u, \mathbf{a}_{u_0})}{\partial u^2}\Big|_{u=u_0}.$$

REMARK 2. (i) We recall the structure of this result: The asymptotic Gaussian distribution is the same as for the stationary approximation. In addition, we have a bias term which comes from the deviation of the true process from the stationary approximation on the segment. In particular, this bias term is zero if the true process is stationary. A simple example is given below. By estimating and minimizing the mean squared error (i.e., by balancing the variance and the bias due to nonstationarity on the segment), we may find an estimator for the optimal segment length.

(ii) If $\mathbb{E}Z_0^4 = 3$, as in the case of normally distributed $Z_t$, then

$$\sqrt{bN}(\hat{\mathbf{a}}_{t_0,N} - \mathbf{a}_{u_0}) + \sqrt{bN}b^2\mu(u_0) \xrightarrow{\mathcal{D}} \mathcal{N}(0, w_2\Sigma(u_0)^{-1}).$$

(iii) It is clear from Propositions 2 and 3 that

$$\mathbb{E}\|\Sigma(u_0)^{-1}\mathcal{B}_{t_0,N}(\mathbf{a}_{u_0})\|_2^2 = b^4\|\mu(u_0)\|_2^2 + O\left(b^6 + \frac{1}{N}\right)$$

and

$$\mathbb{E}\|\Sigma(u_0)^{-1}\nabla\tilde{\mathcal{L}}_N(u_0, \mathbf{a}_{u_0})\|_2^2 = w_2\frac{\operatorname{var}(Z_0^2)}{2bN}\operatorname{trace}(\Sigma(u_0)^{-1}) + o\left(\frac{1}{bN}\right).$$

Therefore, if $b^{13} \ll N^{-1}$, using the above, we conjecture that

$$(37) \qquad \begin{aligned} &\mathbb{E}\|\hat{\mathbf{a}}_{t_0,N} - \mathbf{a}_{u_0}\|_2^2 \\ &= b^4\|\mu(u_0)\|_2^2 + w_2\frac{\operatorname{var}(Z_0^2)}{2bN}\operatorname{trace}(\Sigma(u_0)^{-1}) + o\left(b^4 + \frac{1}{bN}\right). \end{aligned}$$

However, this is very hard to prove. The $b$ which minimizes the conjectured mean square error would be the theoretical optimal bandwidth (i.e., the optimal segment length).

(iv) We illustrate the above results with an example. We first consider the tvARCH(0) process

$$X_{t,N} = \sigma_{t,N}Z_t, \qquad \sigma_{t,N}^2 = a_0\left(\frac{t}{N}\right),$$

which Drees and Stărică [5] have also studied. In this case $\frac{\partial \tilde{X}_t(u)^2}{\partial u} = a_0'(u)Z_t^2$ and under Assumption 2, we have

$$\frac{\partial^2 \nabla \mathcal{L}(u, \mathbf{a}_{u_0})}{\partial u^2}\bigg|_{u=u_0} = -\frac{1}{2}\frac{a_0''(u_0)}{a_0(u_0)^2} \quad \text{and} \quad \Sigma(u_0) = \frac{1}{2a_0(u_0)^2},$$

that is,

$$\mu(u_0) = -\tfrac{1}{2}w(2)a_0''(u_0).$$

This example illustrates well how the bias is linked to the nonstationarity of the process—if the process were stationary, the derivatives of $a_0(\cdot)$ would be zero, causing the bias also to be zero. Conversely, sudden variations in $a_0(\cdot)$ about the time point $u_0$ would be reflected in $a_0''(u_0)$ and manifest as a large $\mu(u_0)$. Straightforward minimization of the first two summands in (37) leads to the optimal bandwidth, which in this case (and for Gaussian $Z_t$) takes the form

$$b_{\text{opt}} = \left(\frac{2w_2}{w(2)^2}\right)^{1/5} N^{-1/5} \left[\frac{a_0(u_0)}{a_0''(u_0)}\right]^{2/5},$$

leading to a large bandwidth if $a_0''(u_0)$ is small and vice versa. Thus, the optimal choice of the bandwidth (of the segment length) depends on the degree of stationarity of the process. For general tvARCH($p$) processes $\mu(u_0)$ is very hard to evaluate. Furthermore, it assumes a very complicated form.

(v) It is of interest to investigate whether the differences in the kernel-QML at each time point are because the true ARCH parameters are time-varying or are simply due to random variation in the estimation method. From a practical point of view, one could evaluate the sum of squared deviations between the kernel-QML estimator at each time point and the global QML estimator. We conjecture that the asymptotic distribution under the null hypothesis of stationarity is a chi-square.

**4. The derivative process.** A key element to the proof of Theorem 3 is the notion of the derivative of the process $\tilde{X}_t(u)^2$ with respect to $u$ and the resulting Taylor expansion for the nonstationary process $X_{t,N}^2$ in terms of stationary processes as given in Corollary 2 below. Since these "derivative processes" are of general interest, we introduce them in this section for general tvARCH($\infty$) processes $X_{t,N}$ as given in (2) and $\tilde{X}_t(u)$ given in (6). We need the following stronger regularity conditions on the parameters.

ASSUMPTION 3. The third derivative of $\{a_j(\cdot)\}$ exists. Furthermore,

$$(38) \qquad \sup_u \left|\frac{\partial^i a_j(u)}{\partial u^i}\right| \le \frac{C}{\ell(j)} \qquad \text{for } i = 1, 2, 3 \text{ and } j = 0, 1, \dots,$$

where $\ell(j)$ is defined as in Assumption 1 and $C$ is a finite constant independent of $i$ and $j$.

THEOREM 4. *Suppose Assumptions 1 and 3 hold and let $\{\tilde{X}_t(u)\}$ be defined as in (6). Then the derivatives $\{\frac{\partial \tilde{X}_t(u)^2}{\partial u}\}$, $\{\frac{\partial^2 \tilde{X}_t(u)^2}{\partial u^2}\}$ and $\{\frac{\partial^3 \tilde{X}_t(u)^2}{\partial u^3}\}$ are almost surely well defined unique stationary stochastic processes for each $u \in (0, 1)$. Furthermore, $\frac{\partial \tilde{X}_t(u)^2}{\partial u}$ is almost surely the unique solution of the stochastic differential equation*

$$(39) \quad \frac{\partial \tilde{X}_t(u)^2}{\partial u} = \left(a_0'(u) + \sum_{j=1}^{\infty} a_j'(u)\tilde{X}_{t-j}(u)^2 + \sum_{j=1}^{\infty} a_j(u)\frac{\partial \tilde{X}_{t-j}(u)^2}{\partial u}\right)Z_t^2,$$

*where $a'_j(u)$ denotes the derivative of $a_j(u)$ with respect to $u$.*

Note that (39) is just the derivative of (6) with $\tilde{\sigma}_t(u)^2$ replaced by $\tilde{X}_t(u)^2/Z_t^2$. An explicit formula for $\frac{\partial \tilde{X}_t(u)^2}{\partial u}$ is given in (51). Similar expressions also hold for the second and third derivatives. For example, if all the derivatives of $a_j(\cdot)$ were zero also the derivative process would be zero [in this case $X_{t,N}^2$ would be stationary and $\tilde{X}_t(u)^2 = X_{t,N}^2$ for all $u$].

An important consequence of Theorem 4 is that it allows us to make a Taylor expansion of $\tilde{X}_t(u)^2$ about $u_0$ (rigorously proved in Section 5), to give

$$
(40) \quad \begin{aligned}
\tilde{X}_t(u)^2 &= \tilde{X}_t(u_0)^2 + (u - u_0)\frac{\partial \tilde{X}_t(u)^2}{\partial u}\Big|_{u=u_0} + \frac{1}{2}(u - u_0)^2\frac{\partial^2 \tilde{X}_t(u)^2}{\partial u^2}\Big|_{u=u_0} \\
&\quad + O_p((u - u_0)^3).
\end{aligned}
$$

An interesting feature of the Taylor expansion in (40) is that it does not depend on the existence of moments of $\tilde{X}_t(u)^2$, unlike other types of series expansions. Instead the expansion depends on the smoothness of the parameters $a_j(\cdot)$.

The approximation in (7), where $X_{t,N}^2 = \tilde{X}_t(\frac{t}{N})^2 + O_p(1/N)$, and the Taylor expansion in (40) lead to the corollary below.

COROLLARY 2. *Suppose $\{X_{t,N}\}$ is a tvARCH process which satisfies Assumptions 1 and 3 and let $\tilde{X}_t(u)$ be defined as in (6). Then for any $u_0 \in (0, 1]$, we have*

$$
(41) \quad \begin{aligned}
X_{t,N}^2 &= \tilde{X}_t(u_0)^2 + \left(\frac{t}{N} - u_0\right)\frac{\partial \tilde{X}_t(u)^2}{\partial u}\Big|_{u=u_0} \\
&\quad + \frac{1}{2}\left(\frac{t}{N} - u_0\right)^2\frac{\partial^2 \tilde{X}_t(u)^2}{\partial u^2}\Big|_{u=u_0} + O_p\left(\left(\frac{t}{N} - u_0\right)^3 + \frac{1}{N}\right).
\end{aligned}
$$

The nice feature of the result of Corollary 2 is that it gives a Taylor expansion of the nonstationary process $X_{t,N}^2$ around $\tilde{X}_t(u_0)^2$ in terms of stationary processes. This is particularly nice since it allows use of well-known results for stationary processes (such as the ergodic theorem) in describing properties of $X_{t,N}$. A similar result also holds for higher-order expansions with higher-order derivatives. However, in this paper only a second-order expansion is needed.

As an example, we now use (41) to derive a tighter bound for the remainder $R_N$ in (8). The effect is similar as in nonparametric regression: Due to the anti-symmetry of the kernel weights, the expectation of the first term falls out. By using (41), we have, for $|t_0/N - u_0| < 1/N$,

$$R_N = \frac{1}{2M+1} \sum_{k=-M}^{M} \frac{k}{N} \frac{\partial \tilde{X}_{t_0+k}(u)^2}{\partial u} \Big|_{u=u_0}$$

$$+ \frac{1}{2M+1} \sum_{k=-M}^{M} \frac{1}{2} \left(\frac{k}{N}\right)^2 \frac{\partial^2 \tilde{X}_{t_0+k}(u)^2}{\partial u^2} \Big|_{u=u_0} + O_p\left(\left(\frac{M}{N}\right)^3 + \frac{1}{N}\right)$$

$$= T_1 + T_2 + O_p\left(\left(\frac{M}{N}\right)^3 + \frac{1}{N}\right).$$

The expectation of $T_1$ is zero. Under the additional condition $\mathbb{E}(Z_0^4)^{1/2} \sum_j Q/\ell(j) \leq (1-\nu)$, $\{\frac{\partial \tilde{X}_t(u)^2}{\partial u}\}$ is a short memory process in which case $\text{var}(T_1) = O(M/N^2)$ (see Lemmas A.7 and A.10). Thus, $T_2$ dominates $T_1$ in probability and we have

$$R_N = \frac{1}{2M+1} \sum_{k=-M}^{M} \frac{1}{2} \left(\frac{k}{N}\right)^2 \frac{\partial^2 \tilde{X}_{t_0+k}(u)^2}{\partial u^2} \Big|_{u=u_0} + O_p\left(\left(\frac{M}{N}\right)^3 + \left(\frac{\sqrt{M}}{N}\right)\right).$$

Note that this is a (stochastic) bias of the approximation in (8).

Theorem 4 and Corollary 2 can easily be generalized to include derivatives of functions of tvARCH processes. By using the chain and product rules, we have the generalization below, which we use to study the quasi-likelihood defined in Section 3.

COROLLARY 3. *Suppose Assumptions 1 and 3 hold, let $\{\tilde{X}_t(u)\}$ be as defined in (6) and $f : \mathbb{R}^d \to \mathbb{R}$, where the first, second and third derivatives of $f$ exist.*

(i) *Then*

$$\frac{\partial f(\tilde{X}_{t_1}(u)^2, \ldots, \tilde{X}_{t_d}(u)^2)}{\partial u} = \sum_{i=1}^{d} \frac{\partial \tilde{X}_{t_i}(u)^2}{\partial u} \frac{\partial f}{\partial \tilde{X}_{t_i}(u)^2},$$

(42)
$$\frac{\partial^2 f(\tilde{X}_{t_1}(u)^2, \ldots, \tilde{X}_{t_d}(u)^2)}{\partial u^2} = \sum_{i=1}^{d} \frac{\partial^2 \tilde{X}_{t_i}(u)^2}{\partial u^2} \frac{\partial f}{\partial \tilde{X}_{t_i}(u_0)^2}$$

$$+ \sum_{i,j=1}^{d} \frac{\partial \tilde{X}_{t_i}(u)^2}{\partial u} \frac{\partial \tilde{X}_{t_j}(u)^2}{\partial u}$$

$$\times \frac{\partial^2 f}{\partial \tilde{X}_{t_i}(u_0)^2 \, \partial \tilde{X}_{t_j}(u_0)^2}.$$

*Furthermore, by using the product and chain rules, similar expressions can be obtained for $\frac{\partial^3 f}{\partial u^3}$.*

(ii) *Suppose $f: \mathbb{R}^d \to \mathbb{R}$ is differentiable with uniformly bounded third derivative. Then we have*

$$
\begin{aligned}
f(X^2_{t+t_1,N}, \ldots, X^2_{t+t_d,N}) = {} & f(\tilde{\mathbf{X}}_t(u)^2) + \left(\frac{t}{N} - u_0\right) \frac{\partial f(\tilde{\mathbf{X}}_t(u)^2)}{\partial u}\Big|_{u=u_0} \\
& + \frac{(t/N - u_0)^2}{2} \frac{\partial^2 f(\tilde{\mathbf{X}}_t(u)^2)}{\partial u^2}\Big|_{u=u_0} \\
& + O_p\left(\left(\frac{t}{N} - u_0\right)^3 + \frac{1}{N}\right),
\end{aligned}
$$

(43)

*where $\tilde{\mathbf{X}}_t(u)^2 := (\tilde{X}_{t+t_1}(u)^2, \ldots, \tilde{X}_{t+t_d}(u)^2)^T$.*

## 5. Volterra expansions of tvARCH processes.

In this section we prove the existence and uniqueness of the process $X_{t,N}$ and of the derivative process from Section 4. This is done by means of Volterra expansions. The methods used here can easily be generalized to include other nonstationary stochastic processes which have as their solution a Volterra expansion. Therefore, the results and methods in this section are of independent interest. A treatise of ordinary Volterra expansions can be found in [14].

Giraitis, Kokoszka and Leipus [9] have shown that a unique solution of $\tilde{X}_t(u)^2$, defined in (6), is almost surely the Volterra series given by

$$
\tilde{X}_t(u)^2 = a_0(u)Z_t^2 + \sum_{k \geq 1} \tilde{m}_t(u, k),
$$

(44)

where

$$
\begin{aligned}
\tilde{m}_t(u, k) &= \sum_{j_1, \ldots, j_k \geq 1} a_0(u) \left(\prod_{r=1}^k a_{j_r}(u)\right) \prod_{r=0}^k Z^2_{t-\sum_{s=1}^r j_s} \\
&= \sum_{j_k < \cdots < j_0 : j_0 = t} \tilde{g}_u(k, j_0, j_1, \ldots, j_k) \prod_{i=0}^k Z^2_{j_i},
\end{aligned}
$$

with

$$
\tilde{g}_u(k, j_0, j_1, \ldots, j_k) = a_0(u) \prod_{i=1}^k a_{(j_{i-1}-j_i)}(u).
$$

We now show a similar result is true for $X^2_{t,N}$. Let $a_j(u) = 0$ for $u < 0$ and $j \geq 0$. A formal expansion of $X^2_{t,N}$, defined in (2), gives

$$
X^2_{t,N} = a_0\left(\frac{t}{N}\right) Z_t^2 + \sum_{k \geq 1} m_{t,N}(k),
$$

(45)

where

$$m_{t,N}(k) = \sum_{j_1,\ldots,j_k \geq 1} a_0 \left( \frac{t - \sum_{s=1}^k j_s}{N} \right) \left( \prod_{r=1}^k a_{j_r} \left( \frac{t - \sum_{s=1}^{r-1} j_s}{N} \right) \right) \left( \prod_{r=0}^k Z_{t-\sum_{s=1}^r j_s}^2 \right)$$

$$= \sum_{j_k < \cdots < j_0 \,:\, j_0 = t} g_{t,N}(k, j_0, j_1, \ldots, j_k) \prod_{i=0}^k Z_{j_i}^2,$$

with

$$g_{t,N}(k, j_0, j_1, \ldots, j_k) = a_0 \left( \frac{j_k}{N} \right) \prod_{i=1}^k a_{(j_{i-1} - j_i)} \left( \frac{j_{i-1}}{N} \right).$$

We stated in Proposition 1 that the tvARCH process has a unique solution. We now prove this result by showing that (45) is the unique solution. The proof in many respects is close to the proof of Theorem 2.1 in [9].

PROOF OF PROPOSITION 1. We first show that (45) is well defined. Since (45) is the sum of positive random variables and the coefficients are also positive, we only need to show that the expectation of (45) is finite. By using (3), (4) and the monotone convergence theorem, a bound for the expectation of (45) is

$$(46) \quad \begin{aligned} \mathbb{E}(X_{t,N}^2) &\leq \sup_u a_0(u) + \sup_u a_0(u) \sum_{k=1}^\infty \sum_{j_k < \cdots < j_0 \,:\, j_0 = t} \prod_{i=1}^k \frac{Q}{\ell(j_i - j_{i-1})} \\ &\leq \sup_u a_0(u) \left[ 1 + \sum_{k=1}^\infty (1 - \nu)^k \right] < \infty. \end{aligned}$$

Furthermore, it is not difficult to see that $X_{t,N}^2$ is a well-defined solution of (2).

To show uniqueness of $X_{t,N}^2$, we must show that any other solution is equal to $X_{t,N}^2$ with probability one. Suppose $Y_{t,N}^2$ is a solution of (2). By recursively applying relation (2) $r$ times to $Y_{t,N}^2$, we have

$$Y_{t,N}^2 = a_0 \left( \frac{t}{N} \right) + \sum_{k=1}^{r-1} m_{t,N}(k) + \sum_{j_r < \cdots < j_0 \,:\, j_0 = t} g_{t,N}(r, j_0, \ldots, j_r) \frac{Y_{j_r,N}^2}{a_0(j_r/N)} \prod_{i=0}^{r-1} Z_{j_i}^2.$$

Thus, the difference between $Y_{t,N}^2$ and $X_{t,N}^2$ is

$$X_{t,N}^2 - Y_{t,N}^2 = A_r - B_r,$$

where

$$A_r = \sum_{k=r}^\infty m_{t,N}(k)$$

and

$$B_r = \sum_{j_r < \cdots < j_0 \, : \, j_0 = t} g_{t,N}(r, j_0, \ldots, j_r) \frac{Y_{j_r,N}^2}{a_0(j_r/N)} \prod_{i=0}^{r-1} Z_{j_i}^2.$$

We now show, for any $\varepsilon > 0$, that $\sum_{r=1}^{\infty} \mathbb{P}(|A_r - B_r| > \varepsilon) < \infty$. By using (3) and (4), we have $\mathbb{E}(A_r) \le C(1 - \nu)^r$. Furthermore, since $Y_{t,N}^2$ is causal, $Y_{j_r,N}^2$ and $\prod_{i=0}^{r-1} Z_{j_i}^2$ are independent (if $i < r$, then $j_i > j_r$). Therefore, $\mathbb{E}(Y_{j_r,N}^2 \prod_{i=0}^{r-1} Z_{j_i}^2) = \mathbb{E}(Y_{j_r,N}^2)$ and we have

$$\mathbb{E}(B_r) = \sum_{j_r < \cdots < j_0 \, : \, j_0 = t} g_{t,N}(r, j_0, \ldots, j_r) \frac{\mathbb{E}(Y_{j_r,N}^2)}{a_0(j_r/N)}$$

$$\le \frac{1}{\inf_u a_0(u)} \sup_{t,N} \mathbb{E}(Y_{t,N}^2)(1 - \nu)^r.$$

Now by using the Markov inequality, we have $\mathbb{P}(A_r > \varepsilon) \le C_1(1 - \nu)^r/\varepsilon$ and $\mathbb{P}(B_r > \varepsilon) \le C_1(1 - \nu)^r/\varepsilon$ for some constant $C_1$. Therefore, $\mathbb{P}(|A_r - B_r| > \varepsilon) \le C_2(1 - \nu)^r/\varepsilon$. Thus, $\sum_{r=1}^{\infty} \mathbb{P}(|A_r - B_r| > \varepsilon) < \infty$ and by the Borel–Cantelli lemma, the event $\{|A_r - B_r| > \varepsilon\}$ can occur only finitely often with probability one. Since this is true for all $\varepsilon > 0$, we have $Y_{t,N} \overset{\text{a.s.}}{=} X_{t,N}$ and therefore the required result.  $\square$

REMARK 3.  It is worth noting that $m_{t,N}(k)$ can be obtained by using the recursion

$$m_{t,N}(k) = Z_t^2 \sum_{j \ge 1} a_j\left(\frac{t}{N}\right) m_{t-j,N}(k-1) \qquad \text{for } k \ge 2,$$

with the initial condition

$$m_{t,N}(1) = Z_t^2 \sum_{j \ge 1} a_j\left(\frac{t}{N}\right) Z_{t-j}^2.$$

Our object now is to prove Theorem 1, that is, to bound the difference between $X_{t,N}^2$ and $\tilde{X}_t(u_0)^2$. More precisely, we will prove under Assumption 1 that

$$(47) \qquad |X_{t,N}^2 - \tilde{X}_t(u_0)^2| \le K\left(\left|\frac{t}{N} - u_0\right| + \frac{1}{N}\right) U_t,$$

where

$$(48) \qquad U_t = Z_t^2 + \sum_{k=1}^{\infty} Q^{k-1} \sum_{j_k < \cdots < j_0 \, : \, j_0 = t} \frac{k|j_0 - j_k|}{\prod_{i=1}^{k} \ell(j_{i-1} - j_i)} \prod_{i=0}^{k} Z_{j_i}^2$$

is a stationary ergodic positive process with finite expectation.

PROOF OF THEOREM 1. To prove (47), we use the triangle inequality to get

$$\left| X_{t,N}^2 - \tilde{X}_t(u_0)^2 \right| \leq \left| X_{t,N}^2 - \tilde{X}_t\left(\frac{t}{N}\right)^2 \right| + \left| \tilde{X}_t\left(\frac{t}{N}\right)^2 - \tilde{X}_t(u_0)^2 \right|$$

and consider bounding $|X_{t,N}^2 - \tilde{X}_t(\frac{t}{N})^2|$ and $|\tilde{X}_t(\frac{t}{N})^2 - \tilde{X}_t(u_0)^2|$ separately. By using (44) and (45), we have

$$
\begin{aligned}
(49) \quad & \left| X_{t,N}^2 - \tilde{X}_t\left(\frac{t}{N}\right)^2 \right| \\
& \leq \sum_{k\geq 1} \sum_{j_k < \cdots < j_0 \,:\, j_0=t} |g_{t,N}(k, j_0, j_1, \ldots, j_k) \\
& \hspace{3cm} - \tilde{g}_{t/N}(k, j_0, j_1, \ldots, j_k)| \prod_{i=0}^{k} Z_{j_i}^2.
\end{aligned}
$$

We notice $j_0/N = t/N$. By successively replacing $a_{(j_{i-1}-j_i)}(\frac{j_{i-1}}{N})$ by $a_{(j_{i-1}-j_i)}(\frac{j_0}{N})$, by using (3), the Lipschitz continuity of the parameters in (5) and that $(j_0 - j_k) \geq (j_0 - j_i)$ (for $i \leq k$) and $(j_0 - j_k) = \sum_{i=1}^{k}(j_{i-1} - j_i)$, we have

$$
\begin{aligned}
(50) \quad & |g_{t,N}(k, j_0, j_1, \ldots, j_k) - \tilde{g}_{t/N}(k, j_0, j_1, \ldots, j_k)| \\
& \leq K Q^{k-1} \frac{k|j_0 - j_k|}{N \prod_{i=1}^{k} \ell(j_i - j_{i-1})},
\end{aligned}
$$

where $K$ is a finite constant. Therefore, by using (49) and (50), we have

$$\left| X_{t,N}^2 - \tilde{X}_t\left(\frac{t}{N}\right)^2 \right| \leq K \frac{1}{N} U_t.$$

Now we bound $|\tilde{X}_t(\frac{t}{N}) - \tilde{X}_t(u_0)^2|$. By using (44), we have

$$
\begin{aligned}
& \left| \tilde{X}_t(u_0)^2 - \tilde{X}_t\left(\frac{t}{N}\right)^2 \right| \\
& \leq \left| a_0(u_0) - a_0\left(\frac{t}{N}\right) \right| Z_t^2 \\
& \quad + \sum_{k\geq 1} \sum_{j_k < \cdots < j_0 \,:\, j_0=t} |\tilde{g}_{u_0}(k, j_0, j_1, \ldots, j_k) - \tilde{g}_{t/N}(k, j_0, j_1, \ldots, j_k)| \prod_{i=0}^{k} Z_{j_i}^2.
\end{aligned}
$$

By using similar methods to those given above, we have

$$\left| \tilde{X}_t(u_0)^2 - \tilde{X}_t\left(\frac{t}{N}\right)^2 \right| \leq K \left| \frac{t}{N} - u_0 \right| U_t.$$

Therefore, we have shown (47).

We now show that $U_t$ is a well-defined stochastic process. Since $U_t$ is the sum of positive random variables, we only need to show that $\mathbb{E}(U_t) < \infty$. Taking the expectation of $U_t$, using (4) and the independence of $\{Z_t^2\}$, we have

$$\mathbb{E}(U_t) = 1 + \sum_{k=1}^{\infty} \sum_{j_k < \cdots < j_0 \,:\, j_0 = t} kQ^{k-1} \frac{|j_0 - j_k|}{\prod_{i=1}^{k} \ell(j_i - j_{i-1})}$$

$$\leq 1 + L \sum_{k=1}^{\infty} k^2 (1-\nu)^{k-1} < \infty,$$

where $L = \sum_{j=1}^{\infty} j/\ell(j)$ [$L$ is finite by definition of $\ell(j)$]. Thus, $\{U_t\}$ is a well-defined process with finite mean. By using Stout [17], Theorem 3.5.8, we can show that $\{U_t\}$ is an ergodic process. Hence, we have the result.  □

We now prove Theorem 4 on the existence of the derivatives of $\tilde{X}_t(u)^2$ with respect to $u$. We will show that this is given by sums of the derivatives of the $\tilde{m}_t(u, k)$ terms in (44), that is,

$$
\begin{aligned}
(51) \quad \frac{\partial \tilde{X}_t(u)^2}{\partial u} &= a_0'(u) Z_t^2 + a_0'(u) \sum_{k \geq 1} \sum_{j_1, \ldots, j_k \geq 1} \left( \prod_{r=1}^{k} a_{j_r}(u) \right) \prod_{r=0}^{k} Z_{t-\sum_{s=1}^{r} j_s}^2 \\
&\quad + a_0(u) \sum_{k \geq 1} \sum_{n=1}^{k} \sum_{j_1, \ldots, j_k \geq 1} a_{j_n}'(u) \left( \prod_{r=1, r \neq n}^{k} a_{j_r}(u) \right) \prod_{r=0}^{k} Z_{t-\sum_{s=1}^{r} j_s}^2.
\end{aligned}
$$

This leads to the Taylor expansions of $\tilde{X}_t(u)^2$ [as stated in (40)] and finally to the Taylor-type representation of $X_{t,N}^2$ stated in Corollary 2. The latter two results are proved below. Throughout the rest of the section $X_{t,N}^2(\omega), \tilde{X}_t(u,\omega)^2$, and so on, denote a specific realization of $X_{t,N}^2, \tilde{X}_t(u)^2$.

PROOF OF THEOREM 4. From (44), we know that $\tilde{X}_t(u)^2$ has almost surely a Volterra series expansion, given by (44), as its unique solution. Therefore, there exists a subset $\mathcal{N}_1(u)$ of the event space where $\mathbb{P}(\mathcal{N}_1(u)^c) = 1$ and

$$
\begin{aligned}
(52) \quad \tilde{X}_t(u,\omega)^2 &= a_0(u) Z_t(\omega)^2 \\
&\quad + a_0(u) \sum_{k \geq 1} \sum_{j_1, \ldots, j_k \geq 1} \left( \prod_{r=1}^{k} a_{j_r}(u) \right) \prod_{r=0}^{k} Z_{t-\sum_{s=1}^{r} j_s}(\omega)^2
\end{aligned}
$$

$\forall \omega \in \mathcal{N}_1(u)^c$. Furthermore, since the random process $\{U_t\}$, defined in (48), is well defined (see Theorem 1), there exists a set $\mathcal{N}_2$ with $\mathbb{P}(\mathcal{N}_2^c) = 1$ and $U_t(\omega)$ finite, for all $\omega \in \mathcal{N}_2^c$. For $\omega \in \mathcal{N}_3(u)^c = \mathcal{N}_1(u)^c \cap \mathcal{N}_2^c$, we consider realizations of the right-hand side of (51) and

$$V_t = \sup_u |a_0'(u)| Z_t^2 + \sup_u |a_0'(u)| \sum_{k \geq 1} \sum_{j_1, \ldots, j_k \geq 1} \left( \prod_{r=1}^k \sup_u a_{j_r}(u) \right) \prod_{r=0}^k Z_{t - \sum_{s=1}^r j_s}^2$$

$$+ \sup_u a_0(u) \sum_{k \geq 1} \sum_{n=1}^k \sum_{j_1, \ldots, j_k \geq 1} \sup_u |a_{j_n}'(u)| \left( \prod_{r=1, r \neq n}^k \sup_u a_{j_r}(u) \right) \prod_{r=0}^k Z_{t - \sum_{s=1}^r j_s}^2.$$

We will now use the following result: Suppose $f(x) = \sum_{j=1}^\infty g_j(x)$ for $x \in [0,1]$, where $f$ is a deterministic function. It is well known if $\sum_{j=1}^\infty g_j'(x)$ is uniformly convergent [which is true if $\sum_{j=1}^\infty \sup_x |g_j'(x)| < \infty$], the derivatives are finite and $\sum_{j=1}^\infty g_j(x)$ converges at least at one point, then $f'(x) = \sum_{j=1}^\infty g_j'(x)$. We now use this result to show that the derivative of $\tilde{X}_t(u)^2$ is well defined. Suppose $\omega \in \mathcal{N}_3(u)^c$. Then by using (52) and (48), we have

$$\tilde{X}_t(u, \omega)^2 \leq \max(1, Q) \sup_u a_0(u) U_t(\omega) < \infty,$$

where the summands in $U_t(\omega)$ are absolutely and uniformly summable. Furthermore, under Assumption 3 and (4), we have, for all $\omega \in \mathcal{N}_3(u)^c$,

$$V_t(\omega) \leq \left( \sup_u |a_0'(u)| + \frac{C}{Q} \sup_u |a_0(u)| \right) U_t(\omega) < \infty.$$

Therefore, $\frac{\partial \tilde{X}_t(u, \omega)^2}{\partial u} = Y_t(u, \omega)$ is almost surely given by (51). By using [17], Theorem 3.5.8, it is clear $\{ \frac{\partial \tilde{X}_t(u)^2}{\partial u} \}$ is an ergodic process.

To show that (51) is the unique solution of (39), we can use the same method as given in the proof of Theorem 1. We omit the details here.

We can use the same method as described above to show that $\{ \frac{\partial^2 \tilde{X}_t(u)^2}{\partial u^2} \}$ and $\{ \frac{\partial^3 \tilde{X}_t(u)^2}{\partial u^3} \}$ are uniquely well-defined ergodic processes. Again, we omit the details.  □

At this point it is easy to derive some moment conditions on $\tilde{X}_t(u)$ and its derivatives.

LEMMA 2.  Suppose $\{X_{t,N} : t = 1, \ldots, N\}$ is a tvARCH($\infty$) process which satisfies Assumptions 1 and 3 and, in addition,

$$\mathbb{E}(Z_0^{2r})^{1/r} \sum_j \frac{Q}{\ell(j)} < (1 - \nu)$$

for $r \geq 1$. Then $\mathbb{E}|U_t|^r < \infty$, $\mathbb{E}|\sup_u \tilde{X}_t(u)^2|^r < \infty$, $\mathbb{E}|\sup_u \frac{\partial \tilde{X}_t(u)^2}{\partial u}|^r < \infty$ and $\mathbb{E}|\sup_u \frac{\partial^2 \tilde{X}_t(u)^2}{\partial u^2}|^r < \infty$ uniformly in $t$.

PROOF. The result follows by applying the Minkowski inequality to (48), (44), (51) and the corresponding formula for the second derivative and using arguments similar to those in (46). We omit the details. □

PROOF OF (40) AND COROLLARY 2. We first prove (40). For $\omega \in \mathcal{N}_3(u)^c \cap \mathcal{N}_3(u_0)^c$, with $\mathcal{N}_3$ as defined above, the Volterra series expansions (44) give solutions of (6) for $\tilde{X}_t(u)^2$ and $\tilde{X}_t(u_0)^2$. The relation (40) now follows from an ordinary Taylor expansion of $\tilde{X}_t(u, \omega)^2$ about $u_0$, noting that $\mathbb{E}(\frac{\partial^3 \tilde{X}_t(u)^3}{\partial u^3}\rfloor_{u=\tilde{U}}) < \infty$ for an arbitrary random variable $\tilde{U}$.

By using (40) and

$$\left| X_{t,N}^2 - \tilde{X}_t\left(\frac{t}{N}\right)^2 \right| \leq \frac{K}{N} U_t,$$

we obtain Corollary 2. □

## 6. Concluding remarks.
We have studied the class of nonstationary ARCH($\infty$) processes with time-varying coefficients. We have shown that, about a given time point, the process can be approximated by a stationary process. Moreover, this approximation has facilitated the Taylor expansion of the tvARCH process in terms of stationary processes. It is worth mentioning that the existence of the derivatives of the coefficients determines the existence of the derivatives of the process and the subsequent Taylor expansion (and not the existence of the moments). The definition of the derivative process and the Taylor expansion is not restricted to tvARCH($\infty$) processes, and with simple modifications can also be applied to other nonstationary processes.

To estimate the time-varying parameters of a time-varying ARCH($p$) ($p < \infty$) process, we have used a weighted quasi-likelihood on a segment. Investigation of the asymptotic properties of the estimator showed an extra bias due to nonstationarity on the segment. This expression can be used to find an adaptive choice of the segment length (by minimizing, e.g., the mean squared error and estimating the second derivative). The relevance of this model for (say) financial data needs further investigation. We conjecture that, by using tvARCH models, the often discussed long range dependence of the squared log returns can be reduced drastically and even disappear completely (there has been some discussion that the long range dependence of the squares is in truth only due to some nonstationarity in the data; see [12]). Furthermore, we conjecture that, for example, the empirical kurtosis of financial log returns is much smaller with a time-varying model than with a classical ARCH model.

Typically for stationary ARCH($p$) processes, the existence of $\mathbb{E}(Z_0^4)$ is assumed in order to show asymptotic normality of the quasi-likelihood estimator. A drawback of our approach is that the expression of the bias given in

(36) holds only under the assumption $[\mathbb{E}(Z_0^{12})]^{1/6} \sum_{j=1}^{p} \frac{Q}{\ell(j)} \leq (1-\nu)$, that is, under the existence of the 12th moment. However, if we assume the weaker condition $\mathbb{E}(Z_0^{4+\delta}) < \infty$, then the segment quasi-likelihood estimator still has asymptotically a normal distribution, but the explicit form of the bias cannot be evaluated (see also Remark 1).

We mention that, unlike the case of stationary GARCH$(p,q)$ models, the time-varying GARCH model is not included in the tvARCH$(\infty)$ class. The investigation of time-varying GARCH$(p,q)$ models is a topic of future research. However, unlike tvGARCH models, the squares of certain tvARCH$(\infty)$ models have "near" long memory behavior (cf. [9, 11]). This is one justification for studying tvARCH$(\infty)$-models.

An important issue not discussed in this paper are the practical aspects when the model is applied. In particular, identifiability requires investigation since both conditional heteroscedasticity and time varying parameters are suitable to model volatility. Theoretically the model is identifiable and we are convinced this also holds in practice for large data sets. However, it has to be checked whether this leads to satisfactory results for moderate sample sizes. Our idea is that the conditional heteroscedasticity models the short term fluctuations, while the time varying parameters model the longer term changes. Of course this can be achieved by a sufficiently large choice of the bandwidth.

## APPENDIX

In this appendix we establish the results required in the proofs of Section 3.

Many of the results related to the local quasi-likelihood defined at (11) depend on the asymptotic limit of the weighted sum of nonstationary random processes. The general method we use to deal with such sums is to substitute an ergodic process for the nonstationary process, and to study the limit of a weighted sum of an ergodic process. In Appendix A.1 we establish results related to the weighted sums of ergodic processes. These results are used in Appendix A.2, where we study the difference between the nonstationary tvARCH process and the corresponding approximating stationary processes. We then use this result to evaluate the limit of weighted sums of functions of tvARCH$(p)$ processes. In Appendix A.3 we investigate the mixing properties of the likelihood process and in Appendix A.4 the bias of the segment estimate from Section 3.

**A.1. Convergence results for weighted sums of random variables.** In this section we prove ergodic type theorems for weighted sums of ergodic processes. In the lemma below we show an almost sure convergence result and in Lemma A.2 we prove convergence in probability for certain triangular arrays.

LEMMA A.1. *Suppose $\{Y_t\}$ is an ergodic sequence with $\mathbb{E}|Y_t| < \infty$ and $W : [-1/2, 1/2] \to \mathbb{R}$ is a kernel function of bounded variation with $\int_{-1/2}^{1/2} W(x) \, dx = 1$. Then*

$$\sum_{k=-M}^{M} \frac{1}{2M+1} W\left(\frac{k}{2M+1}\right) Y_k \xrightarrow{\text{a.s.}} \mu \qquad as \ M \to \infty,$$

*where $\mu = \mathbb{E}(Y_0)$.*

PROOF. Since

$$(A.1) \qquad \sum_{k=-M}^{M} \frac{1}{2M+1} W\left(\frac{k}{2M+1}\right) \to \int_{-1/2}^{1/2} W(x) \, dx = 1,$$

we can assume without loss of generality that $\mu = 0$. We split the sum into negative and positive suffixed elements, which gives

$$(A.2) \quad \begin{aligned} S_M &= \frac{1}{2M+1} \sum_{k=-M}^{0} W\left(\frac{k}{2M+1}\right) Y_k + \frac{1}{2M+1} \sum_{k=1}^{M} W\left(\frac{k+M}{2M+1}\right) Y_k \\ &= N_M + P_M, \end{aligned}$$

and consider first $P_M$. By using summation by parts, we have, with $\mathcal{S}_k = \sum_{i=1}^{k} Y_i$,

$$\begin{aligned} P_M &= \frac{1}{2M+1} \sum_{k=1}^{M-1} \left[ W\left(\frac{k}{2M+1}\right) - W\left(\frac{k+1}{2M+1}\right) \right] \mathcal{S}_k \\ &\quad + \frac{1}{2M+1} W\left(\frac{M}{2M+1}\right) \mathcal{S}_M. \end{aligned}$$

Since $W$ is of bounded variation, this yields

$$|P_M| \le \frac{K}{2M+1} \sup_{k \le M} |\mathcal{S}_k|$$

with some constant $K$. Now the ergodic theorem implies $\mathcal{S}_k(\omega)/k \to 0$ for almost all $\omega$. It is obvious for these $\omega$ that also $P_M(\omega)$ tends to zero. In the same way we conclude that $N_M \to 0$ a.s., which gives the result. $\quad\square$

For kernel estimates about arbitrary center points, the situation is more difficult since we basically average over triangular arrays of observations. We therefore prove in the following lemma only convergence in probability.

LEMMA A.2. *Suppose $\{Y_t\}$ is an ergodic sequence with $\mathbb{E}|Y_t| < \infty$ and $W : [-1/2, 1/2] \to \mathbb{R}$ is a kernel function of bounded variation with*

$\int_{-1/2}^{1/2} W(x) \, dx = 1$. *Then*

$$(A.3) \quad \hat{\mu}_N(u_0) := \sum_{k=p+1}^{N} \frac{1}{bN} W\left(\frac{u_0 - k/N}{b}\right) Y_k \xrightarrow{\mathcal{P}} \mu \qquad for \ u_0 \in [0,1],$$

*where* $b \to 0$, $bN \to \infty$ *as* $N \to \infty$, *and* $\mu = \mathbb{E}(Y_0)$.

PROOF. Again we consider only the case $\mu = 0$. Suppose $N \geq N_0$ with $N_0$ such that $u_0 - p/N_0 > b_0/2$ and $u_0 - 1 < -b_0/2$, $b_0 = b(N_0)$ [i.e., the sum in (A.3) is over the whole domain of $W$]. Let $k_0 = k_0(N)$ be such that $|u_0 - k_0/N| < 1/N$. Since $\{Y_k\}$ is stationary, $\hat{\mu}_N(u_0)$ has the same distribution as

$$\sum_k \frac{1}{bN} W\left(\frac{u_0 - k/N}{b}\right) Y_{k-k_0} = \sum_k \frac{1}{bN} W\left(\frac{u_0 - k_0/N - k/N}{b}\right) Y_k$$

$$= \hat{\mu}_N\left(u_0 - \frac{k_0}{N}\right).$$

Since $W$ is of bounded variation, this is equal to

$$\sum_k \frac{1}{bN} W\left(-\frac{k}{bN}\right) Y_k + R_N$$

with $R_N \leq \frac{K}{bN} \sup_{-bN < k < bN} |Y_k|$. Lemma A.1 implies that the first term converges to zero almost surely [the proof of Lemma A.1 remains the same with $(2M+1)$ replaced by $bN$]. Since $|Y_k| \leq |\mathcal{S}_k| + |\mathcal{S}_{k-1}|$, where $\mathcal{S}_k = \sum_{i=1}^{k} Y_i$, the second term also converges to zero almost surely (as in the proof of Lemma A.1). Therefore,

$$\mathbb{P}(|\hat{\mu}_N(u_0)| \geq \varepsilon) = \mathbb{P}(|\hat{\mu}_N(u_0 - k_0/N)| \geq \varepsilon) \to 0,$$

which gives the result. $\square$

**A.2. Convergence of the local likelihood and its derivatives.** In this section we evaluate the limit of weighted sums of $\{\ell_{t,N}(\boldsymbol{\alpha})\}$, $\{\nabla \ell_{t,N}(\boldsymbol{\alpha})\}$, $\{\nabla^2 \ell_{t,N}(\boldsymbol{\alpha})\}$ and the corresponding stationary approximations. In particular, we prove Lemma 1. Recall the formulas (14)–(16) and the corresponding formulas for $\tilde{\ell}_t(u_0, \boldsymbol{\alpha})$. Let $\kappa = \frac{\rho_2 + 1}{\rho_1}$ and

$$(A.4) \quad \begin{aligned} \Delta_{t,N} &:= \Delta_{t,N}(u_0, U_t, \boldsymbol{\alpha}) \\ &:= \sum_{j=1}^{p} \alpha_j \left\{ \left| \frac{t-j}{N} - u_0 \right| + \frac{1}{N} \right\} U_{t-j}, \end{aligned}$$

with the ergodic process $U_t$ from (48). For a better understanding of the following result, we note that we have, for $|u_0 - \frac{t_0}{N}| < 1/N$, $\mathbb{E}(\Delta_{t_0,N}(u_0, U_t, \boldsymbol{\alpha})) = O(N^{-1})$ and therefore,

$$\Delta_{t_0,N}(u_0, U_t, \boldsymbol{\alpha}) = O_p(N^{-1}),$$

uniformly in $u_0$ and $\boldsymbol{\alpha}$. The same holds for $Z_{t_0}^2 \Delta_{t_0,N}$ since $Z_t^2$ and $\Delta_{t,N}$ are independent with $\mathbb{E}(Z_t^2) = 1$.

In the following lemmas we derive upper bounds for the expressions occurring in (14), (15) and (16) and for the difference between these expressions and the corresponding expressions in $\bar{\ell}_t(u_0, \boldsymbol{\alpha})$ and its derivatives. Assumption 2(iii) immediately yields

$$(A.5) \qquad \frac{\nabla w_{t,N}(\boldsymbol{\alpha})_i}{w_{t,N}(\boldsymbol{\alpha})} \leq \frac{1}{\rho_1}, \qquad \frac{\nabla \tilde{w}_t(u_0, \boldsymbol{\alpha})_i}{\tilde{w}_t(u_0, \boldsymbol{\alpha})} \leq \frac{1}{\rho_1}$$

uniformly in $t$, $N$, $u_0$ and $\boldsymbol{\alpha}$ $(i = 1, \ldots, p+1)$.

LEMMA A.3. *Suppose* $\{X_{t,N} : t = 1, \ldots, N\}$ *is a* tvARCH$(p)$ *process which satisfies Assumption 2(i), (iii). Then*

$$(A.6) \qquad \frac{X_{t,N}^2}{w_{t,N}(\boldsymbol{\alpha})} \leq \kappa Z_t^2 \quad and \quad \frac{\tilde{X}_t(u)^2}{\tilde{w}_t(u, \boldsymbol{\alpha})} \leq \kappa Z_t^2.$$

PROOF. We only prove (A.6) for the tvARCH case; the proof for the stationary case is similar. Since $X_{t,N}^2 = Z_{t,N}^2 \sigma_{t,N}^2$, we have

$$\begin{aligned}
\frac{X_{t,N}^2}{w_{t,N}(\boldsymbol{\alpha})} &= Z_t^2 \left( \frac{\sigma_{t,N}^2}{w_{t,N}(\boldsymbol{\alpha})} \right) \\
&= Z_t^2 \left( \frac{a_0(t/N) + \sum_{j=1}^p a_j(t/N) X_{t-j,N}^2}{\alpha_0 + \sum_{j=1}^p \alpha_j X_{t-j,N}^2} \right) \\
&\leq Z_t^2 \left( \frac{a_0(t/N)}{\alpha_0} + \sum_{j=1}^p \frac{a_j(t/N)}{\alpha_j} \right) \\
&\leq \kappa Z_t^2.
\end{aligned}$$

The last line is true because $\sum_{j=1}^p a_j(t/N) < 1$ and $\alpha_j > \rho_1$ for $j = 0, \ldots, p$. □

LEMMA A.4. *Under the assumptions of Lemma A.3, we have*

$$(A.7) \qquad \frac{X_{t,N}^2}{w_{t,N}(\boldsymbol{\alpha})} = \frac{\tilde{X}_t(u_0)^2}{\tilde{w}_t(u_0, \boldsymbol{\alpha})} + R_{1,N}(u_0, t),$$

$$where \ |R_{1,N}(u_0, t)| \leq \frac{1}{\rho_1} \left( \left| \frac{t}{N} - u_0 \right| \right) U_t + \frac{\kappa}{\rho_1} Z_t^2 \Delta_{t,N}(u_0, U_t, \boldsymbol{\alpha}),$$

$$(A.8) \qquad \frac{\nabla w_{t,N}(\boldsymbol{\alpha})_i}{w_{t,N}(\boldsymbol{\alpha})} = \frac{\nabla \tilde{w}_t(u_0, \boldsymbol{\alpha})_i}{\tilde{w}_t(u_0, \boldsymbol{\alpha})} + R_{2,N}(u_0, t) \qquad (i = 1, \ldots, p+1),$$

$$where \ |R_{2,N}(u_0, t)| \leq \frac{2}{\rho_1^2} \Delta_{t,N}(u_0, U_t, \boldsymbol{\alpha})$$

*and*

(A.9)
$$\log(w_{t,N}(\boldsymbol{\alpha})) = \log(\tilde{w}_t(u_0, \boldsymbol{\alpha})) + R_{3,N}(u_0, t),$$

$$\text{where } |R_{3,N}(u_0, t)| \leq \frac{1}{\rho_1} \Delta_{t,N}(u_0, U_t, \boldsymbol{\alpha}).$$

PROOF. We first prove (A.7). We have

$$|R_{1,N}(u_0, t)| \leq \left| \frac{X_{t,N}^2}{w_{t,N}(\boldsymbol{\alpha})} - \frac{\tilde{X}_t(u_0)^2}{w_{t,N}(\boldsymbol{\alpha})} \right| + \left| \frac{\tilde{X}_t(u_0)^2}{w_{t,N}(\boldsymbol{\alpha})} - \frac{\tilde{X}_t(u_0)^2}{\tilde{w}_t(u_0, \boldsymbol{\alpha})} \right|$$

(A.10)
$$\leq \frac{1}{\rho_1} |X_{t,N}^2 - \tilde{X}_t(u_0)^2|$$

$$+ \frac{\tilde{X}_t(u_0)^2}{\rho_1 \tilde{w}_t(u_0, \boldsymbol{\alpha})} |w_{t,N}(\boldsymbol{\alpha}) - \tilde{w}_t(u_0, \boldsymbol{\alpha})|.$$

From the definitions of $w_{t,N}(\boldsymbol{\alpha})$ and $\tilde{w}_t(u_0, \boldsymbol{\alpha})$ and by using (7), we have

(A.11)
$$|w_{t,N}(\boldsymbol{\alpha}) - \tilde{w}_t(u_0, \boldsymbol{\alpha})| \leq \sum_{j=1}^{p} \alpha_j \left\{ \left| \frac{t-j}{N} - u_0 \right| + \frac{1}{N} \right\} U_{t-j}$$

$$\leq \Delta_{t,N}(u_0, U_t, \boldsymbol{\alpha}).$$

Together with (7) and (A.6), this leads to (A.7). Since $\nabla w_{t,N}(\boldsymbol{\alpha})_i = X_{t+1-i,N}^2$ for $i = 2, \ldots, p+1$, the proof of (A.8) is almost the same, so we omit the details. The case $i = 1$ also follows in the same way.

We now prove (A.9). By differentiating $\log(w_{t,N}(\boldsymbol{\alpha}))$ with respect to $X_{t-j,N}^2$ and using the mean value theorem, we have

(A.12)
$$R_{3,N}(u_0, t) = \frac{1}{\alpha_0 + \sum_{j=1}^{p} \alpha_j Y_j} \sum_{j=1}^{p} \alpha_j (X_{t-j,N}^2 - \tilde{X}_{t-j}(u_0)^2),$$

where $(Y_j : i = 1, \ldots, p)$ are positive random variables [since both $X_{t-i,N}^2$ and $\tilde{X}_{t-i}(u_0)^2$ are positive and $Y_j$ lies in between]. Therefore, by using (7), we have

$$|R_{3,N}(u_0, t)| \leq \frac{1}{\rho_1} \left( \sum_{j=1}^{p} \alpha_j \left\{ \left| \frac{t-j}{N} - u_0 \right| + \frac{1}{N} \right\} U_{t-j} \right) \leq \frac{1}{\rho_1} \Delta_{t,N}(u_0, U_t, \boldsymbol{\alpha}),$$

which is the required result. □

COROLLARY A.1. *Under the assumptions of Lemma* A.3, *we have, for* $n \in \mathbb{N}$,

$$\frac{\prod_{r=1}^{n} \nabla w_{t,N}(\boldsymbol{\alpha})_{i_r}}{w_{t,N}(\boldsymbol{\alpha})^n} = \frac{\prod_{r=1}^{n} \nabla \tilde{w}_t(u_0, \boldsymbol{\alpha})_{i_r}}{\tilde{w}_t(u_0, \boldsymbol{\alpha})^n} + R_{4,N}(u_0, t),$$

(A.13)
$$\text{where } |R_{4,N}(u_0, t)| \leq \frac{2n}{\rho_1^{n+1}} \Delta_{t,N}(u_0, U_t, \boldsymbol{\alpha})$$

*and*

$$(A.14) \quad \frac{X_{t,N}^2 \prod_{r=1}^{n-1} \nabla w_{t,N}(\boldsymbol{\alpha})_{i_r}}{w_{t,N}(\boldsymbol{\alpha})^n} = \frac{\tilde{X}_t(u_0)^2 \prod_{r=1}^{n-1} \nabla \tilde{w}_t(u_0, \boldsymbol{\alpha})_{i_r}}{\tilde{w}_t(u_0, \boldsymbol{\alpha})^n} + R_{5,N}(u_0, t),$$

$$\text{where } |R_{5,N}(u_0, t)| \leq \frac{1}{\rho_1^n}\left(\left|\frac{t}{N} - u_0\right| + \frac{1}{N}\right)U_t + \frac{2\kappa n}{\rho_1^n}Z_t^2 \Delta_{t,N}(u_0, U_t, \boldsymbol{\alpha}),$$

*where $0 < i_r \leq p$ for $r = 1, \ldots, n$.*

PROOF. We can prove (A.13) by successively replacing $\nabla w_{t,N}(\boldsymbol{\alpha})_{i_r}/w_{t,N}(\boldsymbol{\alpha})$ by $\nabla \tilde{w}_t(u_0, \boldsymbol{\alpha})_{i_r}/\tilde{w}_t(u_0, \boldsymbol{\alpha})$ for $r = 1, \ldots, n$. Then by using (A.8) and the bound $\alpha_{i_r} > \rho_1$ for all $\boldsymbol{\alpha} \in \Omega$, we have the result. We can prove (A.14) by using a similar method as above together with (A.6) and (A.7). We omit the details here.  □

LEMMA A.5. *Suppose $\{X_{t,N}\}$ is a tvARCH($p$) process which satisfies Assumption 2(i), (iii) and $W$ is a kernel function of bounded variation with $\int_{-1/2}^{1/2} W(x)\,dx = 1$. Then we have*

$$(A.15) \quad \begin{aligned} &\sum_{k=p+1}^{N} \frac{1}{bN}W\left(\frac{t_0 - k}{bN}\right)\frac{\tilde{X}_k(u_0)^2 \prod_{r=1}^{n-1} \nabla \tilde{w}_k(u_0, \boldsymbol{\alpha})_{i_r}}{\tilde{w}_k(u_0, \boldsymbol{\alpha})^n} \\ &\quad \xrightarrow{\mathcal{P}} \mathbb{E}\left(\frac{\tilde{X}_0(u_0)^2 \prod_{r=1}^{n-1} \nabla \tilde{w}_0(u_0, \boldsymbol{\alpha})_{i_r}}{\tilde{w}_0(u_0, \boldsymbol{\alpha})^n}\right), \end{aligned}$$

$$(A.16) \quad \sum_{k=p+1}^{N} \frac{1}{bN}W\left(\frac{t_0 - k}{bN}\right)\frac{\prod_{r=1}^{n} \nabla \tilde{w}_k(u_0, \boldsymbol{\alpha})_{i_r}}{\tilde{w}_k(u_0, \boldsymbol{\alpha})^n} \xrightarrow{\mathcal{P}} \mathbb{E}\left(\frac{\prod_{r=1}^{n} \nabla \tilde{w}_0(u_0, \boldsymbol{\alpha})_{i_r}}{\tilde{w}_0(u_0, \boldsymbol{\alpha})^n}\right)$$

*and*

$$(A.17) \quad \sum_{k=p+1}^{N} \frac{1}{bN}W\left(\frac{t_0 - k}{bN}\right)\log(\tilde{w}_k(u_0, \boldsymbol{\alpha})) \xrightarrow{\mathcal{P}} \mathbb{E}(\log(\tilde{w}_0(u_0, \boldsymbol{\alpha}))).$$

PROOF. Since $\frac{\nabla \tilde{w}_t(u_0, \boldsymbol{\alpha})_{i_r}}{\tilde{w}_t(u_0, \boldsymbol{\alpha})} \leq 1/\rho_1$, we have by using (A.6) that

$$\frac{\tilde{X}_t(u_0)^2 \prod_{r=1}^{n-1} \nabla \tilde{w}_t(u_0, \boldsymbol{\alpha})_{i_r}}{\tilde{w}_t(u_0, \boldsymbol{\alpha})^n} \leq \frac{\kappa Z_t^2}{\rho_1^n}.$$

By using [17], Theorem 3.5.8, the process $\{\frac{\tilde{X}_t(u_0)^2 \prod_{r=1}^{n-1} \nabla \tilde{w}_t(u_0, \boldsymbol{\alpha})_{i_r}}{\tilde{w}_t(u_0, \boldsymbol{\alpha})^n}\}_t$ is ergodic and by using the bound above has finite mean. By applying Lemma A.2, we have verified (A.15).

Since $\log \rho_1 \leq \log \tilde{w}_t(u_0, \boldsymbol{\alpha}) \leq (\log \rho_2 + (1/\rho_1)\sum_{j=1}^{p} \alpha_j \tilde{X}_{t-j}(u_0)^2)$, $\log \tilde{w}_t(u_0, \boldsymbol{\alpha})$ has a finite mean. (A.16) and (A.17) follow similarly.  □

Lemma A.6. *Under the assumptions of Lemma* A.5 *and* $|t_0/N - u_0| < 1/N$ *with* $u_0 \in (0,1)$, *we have, for all* $n \in \mathbb{N}$,

$$
\text{(A.18)} \quad
\begin{aligned}
\sup_{\boldsymbol{\alpha} \in \Omega} \sum_{k=p+1}^{N} \frac{1}{bN} W\left(\frac{t_0-k}{bN}\right) &\left| \frac{X_{k,N}^2 \prod_{r=1}^{n-1} \nabla w_{k,N}(\boldsymbol{\alpha})_{i_r}}{w_{k,N}(\boldsymbol{\alpha})^n} \right. \\
&\left. - \frac{\tilde{X}_k(u_0)^2 \prod_{r=1}^{n-1} \nabla \tilde{w}_k(u_0, \boldsymbol{\alpha})_{i_r}}{\tilde{w}_k(u_0, \boldsymbol{\alpha})^n} \right| = O_p(b),
\end{aligned}
$$

$$
\text{(A.19)} \quad
\begin{aligned}
\sup_{\boldsymbol{\alpha} \in \Omega} \sum_{k=p+1}^{N} \frac{1}{bN} W\left(\frac{t_0-k}{bN}\right) &\left| \frac{\prod_{r=1}^{n} \nabla w_{k,N}(\boldsymbol{\alpha})_{i_r}}{w_{k,N}(\boldsymbol{\alpha})^n} \right. \\
&\left. - \frac{\prod_{r=1}^{n} \nabla \tilde{w}_k(u_0, \boldsymbol{\alpha})_{i_r}}{\tilde{w}_k(u_0, \boldsymbol{\alpha})^n} \right| = O_p(b)
\end{aligned}
$$

*and*

$$
\text{(A.20)} \quad \sup_{\boldsymbol{\alpha} \in \Omega} \sum_{k=p+1}^{N} \frac{1}{bN} W\left(\frac{t_0-k}{bN}\right) |\log(w_{k,N}(\boldsymbol{\alpha})) - \log(\tilde{w}_k(u_0, \boldsymbol{\alpha}))| = O_p(b).
$$

Proof.    Let

$$
\begin{aligned}
R_N = \sup_{\boldsymbol{\alpha} \in \Omega} \sum_{k=p+1}^{N} \frac{1}{bN} W\left(\frac{t_0-k}{bN}\right) &\left| \frac{X_{k,N}^2 \prod_{r=1}^{n-1} \nabla w_{k,N}(\boldsymbol{\alpha})_{i_r}}{w_{k,N}(\boldsymbol{\alpha})^n} \right. \\
&\left. - \frac{\tilde{X}_k(u_0)^2 \prod_{r=1}^{n-1} \nabla \tilde{w}_k(u_0, \boldsymbol{\alpha})_{i_r}}{\tilde{w}_k(u_0, \boldsymbol{\alpha})^n} \right|.
\end{aligned}
$$

We note first that if $\boldsymbol{\alpha} \in \Omega$, then $\alpha_i \leq \max(1, \rho_2)$, where $\alpha_i$ is the $i$th element of the $(p+1)$-dimensional vector $\boldsymbol{\alpha}$. By using (A.14) and $|\frac{k-j}{N} - u_0| \leq |\frac{k-j}{N} - p|$ when $k$ lies in the support of $W(\frac{t_0-k}{bN})$, we have the bound

$$
|R_N| \leq C\left(b + \frac{p+1}{N}\right) \sum_{k=p+1}^{N} \frac{1}{bN} W\left(\frac{t_0-k}{bN}\right) V_k = C\left(b + \frac{p+1}{N}\right) L_N,
$$

$$
\text{where } V_k = \left\{ U_k + Z_k^2 \sum_{j=1}^{p} U_{k-j} \right\}
$$

and $C$ is a finite constant. Since $U_{t-j}$ (by Theorem 1) and $Z_t^2$ have finite mean and are independent when $j \geq 1$, $\{V_t\}$ has a finite mean. Therefore, by using (A.3), we have that $L_N \xrightarrow{\mathcal{P}} \mathbb{E}(V_0)$ and $|R_N| = O_p(b)$. Thus, we have proved (A.18).

By using (A.13) and (A.9), we can obtain (A.19) and (A.20) similarly. $\square$

Proof of Lemma 1.    We first show (22). To prove uniform convergence, it is sufficient to show both pointwise convergence and equicontinuity in probability of $\tilde{\mathcal{L}}_N(u_0, \boldsymbol{\alpha})$ (since $\Omega$ is compact). By using (A.15) and (A.17), for every $\boldsymbol{\alpha} \in \Omega$, we have

$$\tilde{\mathcal{L}}_N(u_0, \boldsymbol{\alpha}) = \frac{1}{2} \sum_{k=p+1}^{N} \frac{1}{bN} W\left(\frac{t_0 - k}{bN}\right)$$

$$\times \left( \log(\tilde{w}_k(u_0, \boldsymbol{\alpha})) + \frac{\tilde{X}_k(u_0)^2}{\tilde{w}_k(u_0, \boldsymbol{\alpha})} \right)$$

$$\xrightarrow{\mathcal{P}} \mathcal{L}(u_0, \boldsymbol{\alpha}),$$

where $b \to 0$, $bN \to \infty$ as $N \to \infty$. We now show equicontinuity in probability of $\tilde{\mathcal{L}}_N(u_0, \boldsymbol{\alpha})$. By the mean value theorem, for every $\boldsymbol{\alpha}_1, \boldsymbol{\alpha}_2 \in \Omega$, there exists an $\bar{\boldsymbol{\alpha}} \in \Omega$ such that

$$\frac{|\tilde{\mathcal{L}}_N(u_0, \boldsymbol{\alpha}_1) - \tilde{\mathcal{L}}_N(u_0, \boldsymbol{\alpha}_2)|^2}{\|\boldsymbol{\alpha}_1 - \boldsymbol{\alpha}_2\|_2^2}$$

$$\leq \|\nabla \tilde{\mathcal{L}}_N(u_0, \bar{\boldsymbol{\alpha}})\|_2^2$$

$$\leq \frac{1}{2} \sum_{k=p+1}^{N} \frac{1}{bN} W\left(\frac{t_0 - k}{bN}\right) \left\| \left( \frac{\nabla \tilde{w}_k(u_0, \bar{\boldsymbol{\alpha}})}{\tilde{w}_k(u_0, \bar{\boldsymbol{\alpha}})} - \frac{\tilde{X}_k(u_0)^2 \nabla \tilde{w}_k(u_0, \bar{\boldsymbol{\alpha}})}{\tilde{w}_k(u_0, \bar{\boldsymbol{\alpha}})^2} \right) \right\|_2^2.$$

By using (A.5), we have

$$\|\nabla \tilde{\mathcal{L}}_N(u_0, \bar{\boldsymbol{\alpha}})\|_2^2 \leq \sum_{k=p+1}^{N} \frac{1}{bN} W\left(\frac{t_0 - k}{bN}\right) \frac{1}{2\rho_1} \left( 1 + \frac{\tilde{X}_k(u_0)^2}{\rho_1} \right)$$

$$\xrightarrow{\mathcal{P}} \frac{1}{2\rho_1} \mathbb{E}\left( 1 + \frac{\tilde{X}_k(u_0)^2}{\rho_1} \right) < \infty.$$

Therefore, we have that $\tilde{\mathcal{L}}_N(u_0, \cdot)$ is equicontinuous in probability. Now by pointwise convergence of $\tilde{\mathcal{L}}_N(u_0, \boldsymbol{\alpha})$, equicontinuity of $\tilde{\mathcal{L}}_N(u_0, \boldsymbol{\alpha})$ and the compactness of $\Omega$, we have uniform convergence of the kernel quasi-likelihood.

By using (A.18) and (A.20), it is straightforward to verify (23). (24) and (25) can be proved by using the same method as above.    $\square$

**A.3. Mixing properties of the likelihood process.** We now investigate the mixing properties of $\tilde{X}_t(u)^2$ and later $\{\nabla\tilde{\ell}_t(u, \mathbf{a}_{u_0})\}$ and their derivatives with respect to $u$. Our object is to show that the sums of the absolute values of the covariances of the process $\{\nabla\tilde{\ell}_t(u, \mathbf{a}_{u_0})\}$ and its derivatives are finite under suitable regularity conditions. To achieve this, we use a well-known theorem of Gallant and White [8] which states that $\sum_k |\operatorname{cov}(Y_t, Y_{t+k})| < \infty$ if $\{Y_t\}$ is a $L_2$-Near Epoch Dependent ($L_2$-NED) process of size $-\infty$ on the mixing process $\{X_t\}$ of size $-\infty$ (see Lemma A.10 below). To use this result, we need an appropriate mixing process $\{X_t\}$. To this end, we use a result of Basrak, Davis and Mikosch [1], who have shown that a stationary ARCH($p$) process is $\alpha$-mixing with a geometric rate (thus having size $-\infty$) if Assumption 2(iii), (v) is satisfied. Therefore, the stationary ARCH($p$) process $\{\tilde{X}_t(u)^2\}$, under Assumption 2(v), (vi), is $\alpha$-mixing with size $-\infty$. We will use this fact in the following lemmas, where we will show that both the processes $\{\frac{\partial\nabla\tilde{\ell}_t(u, \mathbf{a}_{u_0})}{\partial u}\}_t$ and $\{\frac{\partial^2\nabla\tilde{\ell}_t(u, \mathbf{a}_{u_0})}{\partial u^2}\}_t$ are $L_2$-NED on $\{\tilde{X}_t(u)^2\}_t$.

Let $\mathcal{F}_{t-m}^{t+m} = \sigma(\tilde{X}_{t-m}(u)^2, \ldots, \tilde{X}_{t+m}(u)^2)$ and $\mathbb{E}_{t-m}^{t+m}(Y) = \mathbb{E}(Y | \mathcal{F}_{t-m}^{t+m})$.

LEMMA A.7. *Suppose* $\{X_{t,N} : t = 1, \ldots, N\}$ *is a* tvARCH($p$) *process which satisfies Assumption* 2(i), (iii)–(v).

(i) *If* $\mathbb{E}(Z_0^4)^{1/2} \sum_j \frac{Q}{\ell(j)} < (1 - \nu)$, *then* $\{\frac{\partial\tilde{X}_{t-i}(u)^2}{\partial u}\}_t$ *and* $\{\frac{\partial^2\tilde{X}_{t-i}(u)^2}{\partial u^2}\}_t$ *are $L_2$-NED of size* $-\infty$ *on* $\{\tilde{X}_t(u)^2\}_t$ $(i = 0, \ldots, p)$.

(ii) *If* $\mathbb{E}(Z_0^8)^{1/4} \sum_j \frac{Q}{\ell(j)} < (1 - \nu)$, *then* $\{\frac{\partial\tilde{X}_{t-i}(u)^2}{\partial u} \frac{\partial\tilde{X}_{t-j}(u)^2}{\partial u}\}_t$ *is $L_2$-NED of size* $-\infty$ *on* $\{\tilde{X}_t(u)^2\}_t$ $(i, j = 0, \ldots, p)$.

*Furthermore,* $\{\tilde{X}_t(u)^2\}$ *is $\alpha$-mixing of size* $-\infty$.

PROOF. That $\{\tilde{X}_t(u)^2\}$ is $\alpha$-mixing of size $-\infty$ follows from [1]. We first prove (i) for $i = 0$:

$$(A.21) \qquad \mathbb{E}\left(\frac{\partial\tilde{X}_t(u)^2}{\partial u} - \mathbb{E}_{t-m}^{t+m}\left(\frac{\partial\tilde{X}_t(u)^2}{\partial u}\right)\right)^2 \leq \alpha_m,$$

where $\alpha_m$ has a geometric rate of decay [thus the derivative process is $L_2$-NED of size $-\infty$ on $\{\tilde{X}_t(u)^2\}$]. Since under the quadratic norm $\mathbb{E}_{t-m}^{t+m}(\frac{\partial\tilde{X}_t(u)^2}{\partial u^2})$ is the best projection of $\frac{\partial\tilde{X}_t(u)^2}{\partial u^2}$ onto the sigma algebra $\mathcal{F}_{t-m}^{t+m}$, then

$$(A.22) \qquad \mathbb{E}\left(\frac{\partial\tilde{X}_t(u)^2}{\partial u} - \mathbb{E}_{t-m}^{t+m}\left(\frac{\partial\tilde{X}_t(u)^2}{\partial u}\right)\right)^2 \leq \mathbb{E}\left(\frac{\partial\tilde{X}_t(u)^2}{\partial u} - S\right)^2$$

for all $S \in \mathcal{F}_{t-m}^{t+m}$. We assume from now on that $m > 2p$. Inspired by (51), we now choose

$$S_t^m = a_0'(u)Z_t^2 + \sum_{k=1}^{m-p}\sum_{r=1}^{k+1}\sum_{\substack{t(m)\leq j_k<\cdots<j_0=t \\ j_{k+1}=j_k}} a_{j_{r-1}-j_r}'(u)\left(\prod_{i=1,i\neq r}^{k+1} a_{j_{i-1}-j_i}(u)\right)\prod_{i=0}^{k} Z_{j_i}^2,$$

where $t(m) = t - m + p$ [the index $j_{k+1}$ is introduced to avoid special treatment of $a_0(u)$]. It is clear that $Z_t^2, \ldots, Z_{t-m+p}^2 \in \mathcal{F}_{t-m}^{t+m}$, therefore $S_t^m \in \mathcal{F}_{t-m}^{t+m}$. It is straightforward to show that the following difference can be partitioned as below:

$$\text{(A.23)} \qquad \frac{\partial \tilde{X}_t(u)^2}{\partial u} - S_t^m = A_m + B_m,$$

where

$$A_m = \sum_{k=m+1-p}^{\infty}\sum_{r=1}^{k+1}\sum_{\substack{j_k<\cdots<j_0=t \\ j_{k+1}=j_k}} a_{j_{r-1}-j_r}'(u)\left(\prod_{i=1,i\neq r}^{k+1} a_{j_{i-1}-j_i}(u)\right)\prod_{i=0}^{k} Z_{j_i}^2$$

and

$$B_m = \sum_{k=1}^{m-p}\sum_{r=1}^{k+1}\sum_{\substack{j_k<\cdots<j_0=t \\ j_k<t(m),j_{k+1}=j_k}} a_{j_{r-1}-j_r}'(u)\left(\prod_{i=1,i\neq r}^{k+1} a_{j_{i-1}-j_i}(u)\right)\prod_{i=0}^{k} Z_{j_i}^2.$$

We have

$$\text{(A.24)} \qquad \left\|\frac{\partial \tilde{X}_t(u)^2}{\partial u} - S_t^m\right\|_2 \leq \|A_m\|_2 + \|B_m\|_2.$$

Our object is to show that the mean square error of (A.24) has a geometric rate of decay, which, by the inequality in (A.22), implies (A.21). We now bound $\|A_m\|_2$ and $\|B_m\|_2$. Under Assumption 2(iv), there exists a $C^*$ such that $\sup_u |a_j'(u)| < C^*Q/\ell(j)$ for $j = 1, \ldots, p$. Therefore, by using the Minkowski inequality, (3) and (4), we have

$$\begin{aligned}
\|A_m\|_2 &\leq \sup_{u_1,u_2}(C^*a_0(u_1) + |a_0'(u_2)|) \\
&\qquad \times \sum_{k=m+1-p}^{\infty}\sum_{r=1}^{k}\sum_{j_k<\cdots<j_0=t}\prod_{i=1}^{k}\frac{Q}{\ell(j_{i-1}-j_i)}\mathbb{E}(Z_0^4)^{(k+1)/2} \\
&\leq \sup_{u_1,u_2}(C^*a_0(u_1) + |a_0'(u_2)|)\sum_{k=m+1-p}^{\infty} k(1-\nu)^k \\
&= K(1-\nu)^{m-p},
\end{aligned}$$

(A.25)

where $K$ is a finite constant. Now we bound $\|B_m\|_2$. Since $a_j(u) = 0$ for $j > p$, all $j_{i-1} - j_i$ have to be $\leq p$ in order to have a nonzero contribution in $B_m$. Since for $j_k < t - m + p$,

$$\sum_{i=1}^{k}(j_{i-1} - j_i) = j_0 - j_k \geq m - p,$$

this can only be true for $k \geq (m-p)/p$. Therefore,

$$|B_m| \leq \sum_{k=[(m-p)/p]}^{m-p} \sum_{r=1}^{k+1} \sum_{\substack{j_k < \cdots < j_0 = t \\ j_{k+1} = j_k}} |a'_{j_{r-1} - j_r}(u)| \left( \prod_{i=1, i \neq r}^{k+1} a_{j_{i-1} - j_i}(u) \right) \prod_{i=0}^{k} Z_{j_i}^2,$$

which gives

$$
\begin{aligned}
\|B_m\|_2 &\leq \sup_{u_1, u_2} (C^* a_0(u_1) + |a'_0(u_2)|) \\
&\quad \times \sum_{k=[(m-p)/p]}^{m-p} \sum_{r=1}^{k} \sum_{j_k < \cdots < j_0 = t} \left( \prod_{i=1}^{k} \frac{Q}{\ell(j_{i-1} - j_i)} \right) \mathbb{E}(Z_0^4)^{(k+1)/2} \\
&\leq \sup_{u_1, u_2} (C^* a_0(u_1) + |a'_0(u_2)|) \sum_{k=[(m-p)/p]}^{\infty} k(1-\nu)^k \\
&= K(1-\nu)^{(m-p)/p},
\end{aligned}
$$
(A.26)

where $K$ is a finite constant. Therefore, by using (A.25) and (A.26), we have

$$
\begin{aligned}
\left\| \frac{\partial \tilde{X}_t(u)^2}{\partial u} - \mathbb{E}_{t-m}^{t+m}\left( \frac{\partial \tilde{X}_t(u)^2}{\partial u} \right) \right\|_2 &\leq \left\| \frac{\partial \tilde{X}_t(u)^2}{\partial u^2} - S_t^m \right\|_2 \\
&\leq 2K((1-\nu)^{1/p})^{m-p},
\end{aligned}
$$
(A.27)

thus giving a geometric rate for (A.21) and the required result.

For $\{\frac{\partial \tilde{X}_{t-i}(u)^2}{\partial u}\}$, the result below follows in the same way by using $S_{t-i}^{m-i}$ instead of $S_t^m$. For $\{\frac{\partial \tilde{X}_{t-i}(u)^2}{\partial u} \frac{\partial \tilde{X}_{t-j}(u)^2}{\partial u}\}$, we use the product $S_{t-i}^{m-i} S_{t-j}^{m-j}$ instead of $S_t^m$. Since

$$
\begin{aligned}
&\left\| \frac{\partial \tilde{X}_{t-i}(u)^2}{\partial u} \frac{\partial \tilde{X}_{t-j}(u)^2}{\partial u} - S_{t-i}^{m-i} S_{t-j}^{m-j} \right\|_2 \\
&\quad \leq \left\| \frac{\partial \tilde{X}_{t-i}(u)^2}{\partial u} \right\|_4 \left\| \frac{\partial \tilde{X}_{t-j}(u)^2}{\partial u} - S_{t-j}^{m-j} \right\|_4 \\
&\qquad + \|S_{t-j}^{m-j}\|_4 \left\| \frac{\partial \tilde{X}_{t-i}(u)^2}{\partial u} - S_{t-i}^{m-i} \right\|_4
\end{aligned}
$$

and $\|A_{m-i}\|_4$ and $\|B_{m-i}\|_4$ also have a geometric rate of decay, we also obtain $L_2$-NED of size $-\infty$ in this case. The $L_2$-NED property for $\frac{\partial^2 \tilde{X}_{t-i}(u)^2}{\partial u^2}$ is proved in a similar way. We omit the details. $\quad \square$

In Lemma A.9 we generalize the above result to derivatives of $\{\nabla \tilde{\ell}_t(u, \mathbf{a}_{u_0})\}$ with respect to $u$, which we use in Corollary A.2. We will also need the lemma below, which gives conditions under which moments of $\frac{\partial^s \nabla \tilde{\ell}_t(u, \mathbf{a}_{u_0})_i}{\partial u^s}$ exist.

LEMMA A.8. *Suppose* $\{X_{t,N} : t = 1, \ldots, N\}$ *is a* tvARCH($p$) *process which satisfies Assumption* 2(i), (iii)–(v) *and, in addition,*

$$(A.28) \qquad (\mathbb{E}(Z_0^{2rs}))^{1/rs} \sum_j \frac{Q}{\ell(j)} \leq (1 - \nu)$$

*for* $r \geq 1$ *and* $s \in \mathbb{N}$. *Then*

$$\mathbb{E} \sup_u \left| \frac{\partial^s \nabla \tilde{\ell}_t(u, \mathbf{a}_{u_0})_i}{\partial u^s} \right|^r < \infty \qquad for \ i = 1, \ldots, p + 1,$$

*and the expectation is uniformly bounded in* $u$.

PROOF. We first consider $\nabla \tilde{\ell}_t(u, \mathbf{a}_{u_0})$. It is worth noting

$$(A.29) \qquad \nabla \tilde{\ell}_t(u, \mathbf{a}_{u_0})_1 = \frac{\partial \tilde{\ell}_t(u, a(u_0))}{\partial a_0(u)} = \frac{1}{w_t(u, \mathbf{a}_{u_0})} - \frac{\tilde{X}_t(u)^2}{w_t(u, \mathbf{a}_{u_0})^2}$$

and

$$(A.30) \qquad \begin{aligned} \nabla \tilde{\ell}_t(u, \mathbf{a}_{u_0})_i &= \frac{\partial \tilde{\ell}_t(u, a(u_0))}{\partial a_{i-1}(u)} \\ &= \frac{\tilde{X}_{t-i+1}(u)^2}{w_t(u, \mathbf{a}_{u_0})} - \frac{\tilde{X}_t(u)^2 \tilde{X}_{t-i+1}(u)^2}{w_t(u, \mathbf{a}_{u_0})^2} \end{aligned}$$

for $i = 2, \ldots, p + 1$. We first prove the result for the case $s = 1$. By using Corollary 3(i), we have

$$(A.31) \qquad \frac{\partial \nabla \tilde{\ell}_t(u, \mathbf{a}_{u_0})_i}{\partial u} = \sum_{j=0}^{p} \frac{\partial \tilde{X}_{t-j}(u)^2}{\partial u} \frac{\partial \nabla \tilde{\ell}_t(u, \mathbf{a}_{u_0})_i}{\partial \tilde{X}_{t-j}(u)^2}.$$

By using (A.29) and (A.30), if $\mathbf{a}_{u_0} \in \Omega$, we have

$$(A.32) \qquad \left| \frac{\partial \nabla \tilde{\ell}_t(u, \mathbf{a}_{u_0})_i}{\partial \tilde{X}_t(u)^2} \right| \leq K \quad \text{and} \quad \left| \frac{\partial \nabla \tilde{\ell}_t(u, \mathbf{a}_{u_0})_i}{\partial \tilde{X}_{t-j}(u)^2} \right| \leq K(1 + Z_t^2)$$

$$for \ j = 1, \ldots, p,$$

where $K$ is a finite constant. By using (A.31), (A.32) and the independence of $Z_t^2$ and $\frac{\partial \tilde{X}_{t-j}(u)}{\partial u}$, we obtain

$$\left| \frac{\partial \nabla \tilde{\ell}_t(u, \mathbf{a}_{u_0})_i}{\partial u} \right| \leq \left| \frac{\partial \tilde{X}_t(u)^2}{\partial u} \right| + K \sum_{j=1}^{p} (1 + Z_t^2) \left| \frac{\partial \tilde{X}_{t-j}(u)^2}{\partial u} \right|,$$

thus,

$$(\text{A.33}) \quad \begin{aligned} \left\| \frac{\partial \nabla \tilde{\ell}_t(u, \mathbf{a}_{u_0})_i}{\partial u} \right\|_r &\leq \left\| \frac{\partial \tilde{X}_t(u)^2}{\partial u} \right\|_r + K \sum_{j=1}^p \left\| (1 + Z_t^2) \frac{\partial \tilde{X}_{t-j}(u)^2}{\partial u} \right\|_r \\ &\leq \left\| \frac{\partial \tilde{X}_t(u)^2}{\partial u} \right\|_r + K \sum_{j=1}^p \| 1 + Z_t^2 \|_r \left\| \frac{\partial \tilde{X}_{t-j}(u)^2}{\partial u} \right\|_r. \end{aligned}$$

(A.28) and Lemma 2 now imply the result.

To prove the similar result for the higher order derivatives ($s > 1$), we use the same method as above. But in this case we require stronger conditions on the moments of $\tilde{X}_t(u)^2$ [see (A.28)]. The proof is straightforward and we omit the details here. □

We now use the result above to show that $\{\frac{\partial \nabla \tilde{\ell}_t(u, \mathbf{a}_{u_0})_i}{\partial u}\}$ is $L_2$-NED on $\{\tilde{X}_t(u)^2\}$.

LEMMA A.9.  *Suppose* $\{X_{t,N}\}$ *is a* tvARCH($p$) *process which satisfies Assumption* 2(i), (iii)–(v).

(i) *If* $(\mathbb{E}(Z_0^4))^{1/2} \sum_j \frac{Q}{\ell(j)} \leq (1 - \nu)$, *then the process*

$$\left\{ \frac{\partial \nabla \tilde{\ell}_t(u, \mathbf{a}_{u_0})_i}{\partial u} \right\}_t$$

*is* $L_2$-NED *of size* $-\infty$ *on* $\{\tilde{X}_t(u)^2\}$.

(ii) *If* $(\mathbb{E}(Z_0^8))^{1/4} \sum_j \frac{Q}{\ell(j)} \leq (1 - \nu)$, *then the process*

$$\left\{ \frac{\partial^2 \nabla \tilde{\ell}_t(u, \mathbf{a}_{u_0})_i}{\partial u^2} \right\}_t$$

*is* $L_2$-NED *of size* $-\infty$ *on* $\{\tilde{X}_t(u)^2\}$.

PROOF.  We first prove (i). Let $m > 2p$. By using (A.32), we have

$$\left\| \frac{\partial \nabla \tilde{\ell}_t(u, \mathbf{a}_{u_0})_i}{\partial u} - \mathbb{E}_{t-m}^{t+m}\left( \frac{\partial \nabla \tilde{\ell}_t(u, \mathbf{a}_{u_0})_i}{\partial u} \right) \right\|_2$$

$$\leq \sum_{j=0}^p \left\| \frac{\partial \nabla \tilde{\ell}_t(u, \mathbf{a}_{u_0})_i}{\partial \tilde{X}_{t-j}(u)^2} \left\{ \frac{\partial \tilde{X}_{t-j}(u)^2}{\partial u} - \mathbb{E}_{t-m}^{t+m}\left( \frac{\partial \tilde{X}_{t-j}(u)^2}{\partial u} \right) \right\} \right\|_2$$

$$\leq K \Bigg[ \left\| \frac{\partial \tilde{X}_t(u)^2}{\partial u} - \mathbb{E}_{t-m}^{t+m}\left( \frac{\partial \tilde{X}_t(u)^2}{\partial u} \right) \right\|_2$$

$$+ \| (1 + Z_t^2) \|_2 \sum_{j=1}^p \left\| \left\{ \frac{\partial \tilde{X}_{t-j}(u)^2}{\partial u} - \mathbb{E}_{t-m}^{t+m}\left( \frac{\partial \tilde{X}_{t-j}(u)^2}{\partial u} \right) \right\} \right\|_2 \Bigg].$$

Lemma A.7 now implies that $\{\frac{\partial \nabla \tilde{\ell}_t(u, \mathbf{a}_{u_0})_i}{\partial u}\}_t$ is $L_2$-NED of size $-\infty$ on $\{\tilde{X}_t(u)^2\}_t$.

The proof for the second derivative process is similar, but requires the stronger moment condition given in (ii). We omit the details of the proof. $\square$

We now state a theorem of Gallant and White [8] which we use in Corollary A.2.

LEMMA A.10.  *Suppose the stationary process $\{Y_t\}$ is $L_2$-NED of size $-\infty$ on $\{X_t\}$, which is an $\alpha$-mixing process of size $-\infty$, and we have $\mathbb{E}(Y_0^{2+\delta}) < \infty$ for some $\delta > 0$. Then*

$$\sum_{s=0}^{\infty} |\operatorname{cov}(Y_t, Y_{t+s})| < \infty.$$

COROLLARY A.2.  *Suppose $\{X_{t,N} : t = 1, \ldots, N\}$ is a tvARCH process which satisfies Assumption 2(i), (ii), (iv) and (v) and, in addition, for some $\delta > 0$:*

(i) *If $(\mathbb{E}(Z_0^{2(2+\delta)}))^{1/(2+\delta)} \sum_j \frac{Q}{\ell(j)} \leq (1 - \nu)$, then we have*

(A.34)    $$\sum_{s=0}^{\infty} \left| \operatorname{cov}\left( \frac{\partial \nabla \tilde{\ell}_t(u, \mathbf{a}_{u_0})_i}{\partial u}, \frac{\partial \nabla \tilde{\ell}_{t+s}(u, \mathbf{a}_{u_0})_i}{\partial u} \right) \right| < \infty.$$

(ii) *If $(\mathbb{E}(Z_0^{2(4+\delta)}))^{1/(4+\delta)} \sum_j \frac{Q}{\ell(j)} \leq (1\nu)$, then we have*

(A.35)    $$\sum_{s=0}^{\infty} \left| \operatorname{cov}\left( \frac{\partial^2 \nabla \tilde{\ell}_t(u, \mathbf{a}_{u_0})_i}{\partial u^2}, \frac{\partial^2 \nabla \tilde{\ell}_{t+s}(u, \mathbf{a}_{u_0})_i}{\partial u^2} \right) \right| < \infty.$$

PROOF.  The condition $(\mathbb{E}(Z_0^{2(2+\delta)}))^{1/(2+\delta)} \sum_j \frac{Q}{\ell(j)} \leq (1 - \nu)$ implies $(\mathbb{E}(Z_0^4))^{1/2} \sum_j \frac{Q}{\ell(j)} \leq (1 - \nu)$ by Hölder's inequality. Now by Lemma A.9(i), we have under this assumption that $\{\frac{\nabla \tilde{\ell}_t(u, \mathbf{a}_{u_0})_i}{\partial u}\}_t$ is $L_2$-NED of size $-\infty$ on the $\alpha$-mixing process $\{\tilde{X}_t(u_0)^2\}_t$. Therefore, all the conditions in Lemma A.10 are satisfied and (i) follows. The proof of (ii) is the same, but the stronger condition given in (ii) is required.  $\square$

**A.4. The bias of the segment quasi-likelihood estimate.**

PROOF OF PROPOSITION 3.  Substituting (34) into (33) gives

$$\nabla \mathcal{B}_{t_0, N}(\mathbf{a}_{u_0}) = \tilde{A}_N\left(\frac{t_0}{N}\right) + \frac{1}{2} \tilde{B}_N\left(\frac{t_0}{N}\right) + \frac{1}{3!} \tilde{C}_N\left(\frac{t_0}{N}\right) + R_N,$$

where

$$\tilde{A}_N\left(\frac{t_0}{N}\right) = \sum_k \frac{1}{bN} W\left(\frac{t_0-k}{bN}\right)\left(\frac{k}{N}-u_0\right)\frac{\partial\nabla\tilde{\ell}_k(u,\mathbf{a}_{u_0})}{\partial u}\Big|_{u=u_0},$$

$$\tilde{B}_N\left(\frac{t_0}{N}\right) = \sum_k \frac{1}{bN} W\left(\frac{t_0-k}{bN}\right)\left(\frac{k}{N}-u_0\right)^2\frac{\partial^2\nabla\tilde{\ell}_k(u,\mathbf{a}_{u_0})}{\partial u^2}\Big|_{u=u_0},$$

$$\tilde{C}_N\left(\frac{t_0}{N}\right) = \sum_k \frac{1}{bN} W\left(\frac{t_0-k}{bN}\right)\left(\frac{k}{N}-u_0\right)^3\frac{\partial^3\nabla\tilde{\ell}_k(u,\mathbf{a}_{u_0})}{\partial u^3}\Big|_{u=\tilde{U}_k}.$$

We now consider the expectation of $\tilde{A}_N(u_0)$. We have

$$\mathbb{E}(\tilde{A}_N(u_0)) = \mathbb{E}\left(\frac{\partial\nabla\tilde{\ell}_k(u,\mathbf{a}_{u_0})}{\partial u}\Big|_{u=u_0}\right)\sum_k \frac{1}{bN} W\left(\frac{t_0-k}{bN}\right)\left(\frac{k}{N}-u_0\right)$$

$$= \mathbb{E}\left(\frac{\partial\nabla\tilde{\ell}_k(u,\mathbf{a}_{u_0})}{\partial u}\Big|_{u=u_0}\right)\int_{-1/2}^{1/2}\frac{1}{b}W\left(\frac{x}{b}\right)x\,dx + O\left(\frac{1}{N}\right)$$

$$= O\left(\frac{1}{N}\right).$$

Furthermore, we have

$$\text{var}(\tilde{A}_N(u_0)) = \frac{1}{(bN)^2}\sum_{k_1}\sum_{k_2} W\left(\frac{t_0-k_1}{bN}\right)$$

$$\times W\left(\frac{t_0-k_2}{bN}\right)\left(\frac{k_1}{N}-u_0\right)\left(\frac{k_2}{N}-u_0\right)$$

$$\times \text{cov}\left(\frac{\partial\nabla\tilde{\ell}_{k_1}(u_0,\mathbf{a}_{u_0})}{\partial u}\Big|_{u=u_0},\frac{\partial\nabla\tilde{\ell}_{k_2}(u,\mathbf{a}_{u_0})}{\partial u}\Big|_{u=u_0}\right)$$

$$\leq \frac{b^2}{(bN)^2}\sum_k W\left(\frac{t_0-k}{bN}\right)\sum_s W\left(\frac{t_0-k-s}{bN}\right)$$

$$\times\left|\text{cov}\left(\frac{\partial\nabla\tilde{\ell}_k(u,\mathbf{a}_{u_0})}{\partial u}\Big|_{u=u_0},\frac{\partial\nabla\tilde{\ell}_{k+s}(u,\mathbf{a}_{u_0})}{\partial u}\Big|_{u=u_0}\right)\right|$$

$$\leq \frac{b^2\|W\|_\infty}{(bN)^2}\sum_k W\left(\frac{t_0-k}{bN}\right)$$

$$\times\sum_s\left|\text{cov}\left(\frac{\partial\nabla\tilde{\ell}_k(u,\mathbf{a}_{u_0})}{\partial u}\Big|_{u=u_0},\frac{\partial\nabla\tilde{\ell}_{k+s}(u,\mathbf{a}_{u_0})}{\partial u}\Big|_{u=u_0}\right)\right|,$$

where $\|W\|_\infty = \sup_x W(x)$. By using Corollary A.2, we have that the sum of the absolute values of the covariances is finite. This gives

$$\mathrm{var}(\tilde{A}_N(u_0)) \leq \frac{b^2 \|W\|_\infty^2}{bN}$$

$$\times \sum_s \left| \mathrm{cov}\left( \frac{\partial \nabla \tilde{\ell}_k(u, \mathbf{a}_{u_0})}{\partial u} \Big|_{u_0}, \frac{\partial \nabla \tilde{\ell}_{k+s}(u, \mathbf{a}_{u_0})}{\partial u} \Big|_{u_0} \right) \right|$$

$$= O\left( \frac{b^2}{bN} \right) = O\left( \frac{1}{N} \right).$$

In the same way we obtain

$$\mathbb{E}\left( \tilde{B}_N\left( \frac{t_0}{N} \right) \right) = b^2 w(2) \frac{\partial^2 \nabla \mathcal{L}(u, \mathbf{a}_{u_0})}{\partial u^2} \Big|_{u=u_0} + O\left( \frac{1}{N} \right)$$

and

$$\mathrm{var}\left( \tilde{B}_N\left( \frac{t_0}{N} \right) \right) = O\left( \frac{b^4}{bN} \right) = O\left( \frac{1}{N} \right).$$

We now evaluate a bound for $\mathbb{E}(\tilde{C}_N(u_0)^2)$, which will help us to bound both $\mathbb{E}(\tilde{C}_N(u_0))$ and $\mathrm{var}(\tilde{C}_N(u_0))$. By using Lemma A.8, we have

$$\mathbb{E}(\tilde{C}_N(u_0)^2) = \frac{1}{(bN)^2} \sum_{k_1} \sum_{k_2} W\left( \frac{t_0 - k_1}{bN} \right) W\left( \frac{t_0 - k_2}{bN} \right) \left( \frac{k_1}{N} - u_0 \right)^3 \left( \frac{k_2}{N} - u_0 \right)^3$$

$$\times \mathbb{E}\left( \frac{\partial^3 \nabla \tilde{\ell}_{k_1}(u, \mathbf{a}_{u_0})}{\partial u^3} \Big|_{u=\tilde{U}_{k_1}} \frac{\partial^3 \nabla \tilde{\ell}_{k_2}(u, \mathbf{a}_{u_0})}{\partial u^3} \Big|_{u=\tilde{U}_{k_2}} \right)$$

$$\leq \frac{b^6}{(bN)^2} \mathbb{E}\left( \sup_u \left( \frac{\partial^3 \nabla \tilde{\ell}_k(u, \mathbf{a}_{u_0})}{\partial u^3} \right)^2 \right)$$

$$\times \sum_{k_1} \sum_{k_2} W\left( \frac{t_0 - k_1}{bN} \right) W\left( \frac{t_0 - k_2}{bN} \right)$$

$$\leq b^6 \|W\|_\infty^2 \mathbb{E}\left( \sup_u \left( \frac{\partial^3 \nabla \tilde{\ell}_k(u, \mathbf{a}_{u_0})}{\partial u^3} \right)^2 \right) = O(b^6),$$

leading to the result.  □

**Acknowledgments.** We are grateful to Angelika Rohde for her contributions to the proofs in Appendix A.1, to Sébastien Van Bellegem for careful reading of the manuscript and an Associate Editor and two anonymous referees for several helpful comments.

# REFERENCES

[1] Basrak, B., Davis, R. A. and Mikosch, T. (2002). Regular variation of GARCH processes. *Stochastic Process. Appl.* **99** 95–115. MR1894253

[2] Berkes, I., Horváth, L. and Kokoszka, P. (2003). GARCH processes: Structure and estimation. *Bernoulli* **9** 201–227. MR1997027

[3] Bollerslev, T. (1986). Generalized autoregressive conditional heteroscedasticity. *J. Econometrics* **31** 307–327. MR0853051

[4] Dahlhaus, R. (1997). Fitting time series models to nonstationary processes. *Ann. Statist.* **25** 1–37. MR1429916

[5] Drees, H. and Stărică, C. (2003). A simple non-stationary model for stock returns. Preprint.

[6] Engle, R. F. (1982). Autoregressive conditional heteroscedasticity with estimates of the variance of United Kingdom inflation. *Econometrica* **50** 987–1007. MR0666121

[7] Feng, Y. (2005). Modelling of scale change, periodicity and conditional heteroskedasticity in returns volatility. Preprint. Available at www.ma.hw.ac.uk/~yuanhua/papers/.

[8] Gallant, A. R. and White, H. (1988). *A Unified Theory of Estimation and Inference for Nonlinear Dynamic Models.* Basil Blackwell, Oxford.

[9] Giraitis, L., Kokoszka, P. and Leipus, R. (2000). Stationary ARCH models: Dependence structure and central limit theorem. *Econometric Theory* **16** 3–22. MR1749017

[10] Hall, P. and Heyde, C. (1980). *Martingale Limit Theory and Its Application.* Academic Press, New York. MR0624435

[11] Koulikov, D. (2003). Long memory ARCH($\infty$) models: Specification and quasi-maximum likelihood estimation. Working Paper 163, Centre for Analytical Finance, Univ. Aarhus. Available at www.cls.dk/caf/wp/wp-163.pdf.

[12] Mikosch, T. and Stărică, C. (2000). Is it really long memory we see in financial returns? In *Extremes and Integrated Risk Management* (P. Embrechts, ed.) 149–168. Risk Books, London.

[13] Nelson, D. (1991). Conditional heteroskedasity in assit returns: A new approach. *Econometrica* **59** 347–370. MR1097532

[14] Priestley, M. B. (1988). *Non-Linear and Non-Stationary Time Series Analysis.* Academic Press, London. MR0991969

[15] Robinson, P. M. (1989). Nonparametric estimation of time-varying parameters. In *Statistical Analysis and Forecasting of Economic Structural Change* (P. Hackl, ed.) 253–264. Springer, Berlin.

[16] Robinson, P. M. (1991). Testing for strong serial correlation and dynamic conditional heteroskedasticity in multiple regression. *J. Econometrics* **47** 67–84. MR1087207

[17] Stout, W. (1974). *Almost Sure Convergence.* Academic Press, New York. MR0455094

Institut für Angewandte Mathematik
Universität Heidelberg
Im Neuenheimer Feld 294
69120 Heidelberg
Germany
E-mail: dahlhaus@statlab.uni-heidelberg.de
suhasini@statlab.uni-heidelberg.de